\documentclass[10pt]{article}
\usepackage{fancyhdr}
\usepackage[english]{babel}

\usepackage{amsmath,amsfonts,amssymb}
\usepackage{mathrsfs}
\usepackage{amsthm}
\usepackage{enumerate}
\usepackage[utf8]{inputenc}
\usepackage{fontenc}
\usepackage{xfrac}

\usepackage{amsbsy,amscd,
	graphicx}
\usepackage{hyperref}

\usepackage{colortbl}
\usepackage{AV_gnl_nota}
\usepackage{packUQSD_v2_man}

\begin{document}
	
	\title{Two level natural selection \\with a quasi-stationarity approach}
	\author{Aurélien Velleret\footnote{Aix-Marseille Université, CNRS, Centrale Marseille, I2M, UMR 7373 13453
			Marseille, France, 
			email~: aurelien.velleret@univ-amu.fr}}
	\maketitle	
	
	\setcounter{eq}{0}

\section*{Abstract}
\setcounter{eq}{0}

In a view for a simple model 
where natural selection at the individual level
is confronted to selection effects at the group level,
we consider some individual-based models
of some large population subdivided 
into a large number of groups.
We then obtain the convergence 
to the law of a stochastic process
with some Feynman-Kac penalization.
To analyze the limiting behavior of this law,
we exploit a recent approach, 
designed for the convergence to quasi-stationary distributions.
We are able to deal with the fixation of the stochastic process
and relate the convergence to equilibrium 
to the one where fixation implies extinction.
We notably establish different regimes of convergence.
Besides the case of an exponential rate 
(the rate being uniform over the initial condition),
critical regimes with convergence in $1/t$ are also to notice.
We finally address the relevance of such limiting behaviors 
to predict the long-time behavior of the individual-based model
and describe more specifically the cases of weak selection.
Consequences in term of evolutive dynamics are also derived,
where such competition is assumed to occur repeatedly
at each de novo mutation.



\section*{Introduction}
\setcounter{eq}{0}


We consider a model
of two alleles competing 
in groups of individuals
without inter-group migration.
This model is derived 
as a new limit of large population 
(both within and between groups)
from the more realistic individual-based model 
presented in \cite{L13},
so as to shed light on the dynamics of the latter.
Different scenarii can be observed
depending on the effect of alleles
on the replication within each group
of individuals carrying it
and on the replication of groups as a whole
(where the groups duplicate identically).
The focus is especially 
on two conflicting behaviors :
either the allele favoring replication at group level
is less favorable to the individuals carrying it
for the competition within groups
(case of an altruistic trait);
or  selective effects at group level favors polymorphism
while those at individual level favors a specific allele.

We study a specific limiting behavior 
when the populations sizes are large
(both for the number of groups
and of individuals within each group).
Contrary to \cite{LM15},
we allow for stochastic neutral fluctuations 
of type frequencies within the groups
and show that it leads in fact
to non-trivial effects of selection.
Notably, the strength of selection between groups 
depends very much 
upon sufficiently high levels of stochasticity
for the dynamics within groups.

This mainly explains the discrepancy
between our results 
and the ones of \cite{LM15},
where only the selective effects are kept to the limit.
Considering the specifications of their model,
we look at the competition between an altruistic allele 
and an egoistic allele. 
Type $D$ individuals (the egoistic "Defectors")
perform always better inside their group
while type $C$ individuals (the altruistic "Cooperators")
enhance the global survival of their group.
In our case, 
the equilibrium with a domination of altruistic groups 
is shown to possibly have a much larger basin of attraction :
contrary to the purely selective case (described in \cite{LM15}),
it may reasonably be the endpoint of the dynamics
even though the altruistic allele
constitutes initially a majority in none of the groups.
By reasonably, 
we mean that the convergence
is likely to be seen in the individual-based models,
as we discuss in Section \ref{TS:sec:Disc}.
The crucial issue for the relevance of this approximation
is the following :
the emergence of pure groups
of a mutant type is unlikely
if the invasion of one group
by  individuals of this mutant type
is too exceptional an event 
to be reasonably expected in the population of $m$ groups,
cf Section \ref{TS:sec:pureSel}.

Of course, 
it also highly depends on the range of parameters involved.
Namely, 
when the level of within group stochasticity is small,
we shall retrieve a sharp transition 
similar to the one obtained in \cite{LM15}.
The basin of attraction of altruistic domination 
is in practice limited to initial conditions
where a non-negligible proportion of groups 
have a large majority of $C$ individuals.
The rest of the initial conditions 
belongs to the basin of attraction of defection,
where $D$ individuals have fixed in every group.
\\

In Section \ref{TS:sec:Disc},
we will provide
a much more detailed presentation 
of  the implication  of our results
to describe the dynamics of the individual-based model.
This discussion is a significant part of this work
and clarifies the implications of our rigorous mathematical analysis
of the asymptotic model.
In this section \ref{TS:sec:Disc}, 
we will also discuss implications of our study
in terms of evolutionary dynamics with well-separated mutations :
the central question is then to quantify the probability
of invasion when one type (the mutant)
starts with only one individual
in an otherwise homogeneous population (of residents).

But first of all, 
we will introduce in much more details
both the individual-based models (IBM) under consideration
and several asymptotic descriptions 
in the limit of large population sizes.
These descriptions takes the form of measure-valued processes
that approximate the empirical distribution in the groups
of the proportion of one type. 
This will be done in Section \ref{TS:sec:IBM}
for the description of the individual-based model
and in Section  \ref{TS:sec:Lim} 
for the definitions of these processes 
and their relation to the IBM.
In Section \ref{TS:sec:MCNE}, 
we focus on a different characterization 
of the process
that will be the main subject of study
for the propositions presented in Section \ref{TS:sec:QSD}.
It is expressed
in terms of the law of a stochastic process
with a specific conditioning.
The propositions which follow
provide an exhaustive description 
of the long-time behavior
of this asymptotic solution,
with a variety of potential behaviors.
Some subsequent results for parameters going either to $0$ or infinity
are also given, so as to introduce the already-mentioned discussion
of Section \ref{TS:sec:Disc}.

The proofs of the convergence result to our asymptotic solution
and of the propositions of Section \ref{TS:sec:QSD}
are deferred respectively to Section \ref{TS:PMLim}
and Section \ref{TS:sec:PfQSD}.
For clarity, the latter are kept in the order of appearance of the Theorems,
while the more classical techniques at use in Section \ref{TS:PMLim}
deserve less attention.

\section{The derivation of the limiting model from the individual-based model}
\setcounter{eq}{0}

\subsection*{The individual-based model}
\label{TS:sec:IBM}
Each group contains $n \in \N$ individuals. 
There are two types of individuals: 
type D individuals (for Defectors)
are selectively advantageous 
at  individual level 
while type C individuals (for Cooperators) 
are favorable to the group to which they belong.
Replication and selection occur concurrently 
at the individual and group level according to the Moran process,
as presented in Figure 1 of \cite{LM15}. 
Type C individuals replicate 
at rate $\bar{\gamma}_I \ge 0$
and type D individuals 
at rate $\bar{\gamma}_I \,(1 + \bar{s}), \bar{s}$ 
with $\bar{s}\ge 0$. 
When an individual gives birth, another individual in the same
group is selected uniformly at random to die. 
To reflect the antagonism at the higher level of selection, 
groups replicate at a rate 
which depends on the proportion of type C individuals they contain.
As a simple case, 
we take this rate 
to be of the form
$\bar{\gamma}_G\times \,[1 + \bar{r}(k/n)]$, 
where $k/n$ is the fraction of individuals in the group 
that are of type C
and 
$\bar{r}(x), x\in [0,1]$ is a non-negative bounded function 
measuring selective advantage at group level.
Similarly as at individual level, 
the number of groups 
is maintained at the value $m$ 
by selecting a group uniformly at random to die 
whenever a group replicates. 
The offspring of groups are assumed to be identical to their parent.
We refer to \cite{L13} for a general presentation 
of the biological motivations for such models.

Let $X^i_t$
be the number of type C individuals in group i at time t. Then
\begin{align*}
		\mu^{m;n}_t 
		:= \divi{m} \sum_{i\le m} \delta_{X^i_t/n}
\end{align*}
is the empirical measure at time $t$ 
--of the proportion of type C by group--
for a given number of groups m and individuals per group $n$. 
$\delta_x(y) = 1$ if $x = y$ and zero otherwise. 
The $X^i_t$ are divided by $n$ so that $\mu^{m;n}_t$  is a
probability measure on $E_n := [0;1/n; ... ; 1].$
For fixed $T > 0$, $\mu^{m;n}_t \in D([0; T ]; \M_1(E_n))$, 
the set of càdlàg processes on $[0; T ]$
taking values in $\M_1(E_n)$, 
where $\M_1(S)$ is the set of probability measures on a set $S$. 
The particle
process being as described above, 
the generator of the Markov process $\mu^{m;n}_t$ takes the form :
\begin{align*}
		(\cL^{m;n} \psi )(v) 
		= \sum_{i,j}
		(\bar{\gamma}_I\, R_I^{i,j} + \bar{\gamma}_G\, R_G^{i,j})(v)
		\times \lp  
		\psi\big [v+ 1/m\, (\delta_{j/n} - \delta_{i/n})\big] - \psi[v]\rp
\end{align*}
where $\psi \in C_b(\M_1([0; 1]))$ is a bounded continuous functions, 
and $v \in \M_1(E_n) \subset \M_1([0; 1]).$ 
The transition rates $(\bar{\gamma}_I\, R_I^{i,j} + \bar{\gamma}_G\,R_G^{i,j})$ are given by
\begin{align*}
		R_I^{i,j} (v) 
		:= \left\{
		\begin{aligned}
			&m\, v(i/n) \, i\, (1 - i/n)\, (1+\bar{s}) 
			&\text{ if } j = i - 1; i < n,
			\\& m\, v(i/n) \, i\, (1 - i/n) 
			&\text{ if } j = i + 1; i > 0
			\\& 0 \text{ otherwise}
		\end{aligned}
		\right.
		\EQn{TS:RI}{R_I}
\end{align*}
\begin{align*}
		\AND \hcm{2}
		R_G^{i,j}  (v) 
		:=  m\, v(i/n) \,v(j/n)\, (1 + \bar{r}[j/n] ).
		\EQn{TS:RG}{R_G}
\end{align*}

$R_I^{i,j}$ represents individual-level events 
while $R_G^{i,j}$ represents group-level events.

\subsection{Different limiting behaviors}
\label{TS:sec:Lim}

In \cite{LM15},
Luo and Mattingly  
consider extensively 
the limit as $n, m \ifty$ 
of the measure valued process ($\mu^{m;n}_t$)
with fixed parameters $\bar{\gamma}_I$, $\bar{\gamma}_G$, $\bar{s}$ and $\bar{r}$.
The limiting measure $\pi_t$ satisfies~:
\begin{align*}
		\partial_t\, \LAg\pi_t\bv f\RAg
		= - \bar{\gamma}_I\, s \LAg \pi_t\bv x(1-x)f'\RAg
		+ \bar{\gamma}_G\, [\LAg \pi_t\bv r\,f\RAg  
		- \LAg \pi_t\bv f\RAg \; \LAg \pi_t\bv r\RAg].
		\EQn{TS:piT}{\pi_t}
\end{align*}
They also proved that,
with the scaling :
$\bar{\gamma}_I = n\, \gamma_I$, 
$\bar{\gamma}_G = m\,\gamma_G$,
$\bar{\gamma}_I\, \bar{s} = s$, 
$\bar{\gamma}_G\, \bar{r} = r$,
as $n, m\ifty$,
the process ($\mu^{m;n}_t$) converges weakly to $\nu_t$,
where $\nu_t$ satisfies the following
martingale problem~:
\begin{align*}
		\EQn{TS:LWF}{\cL_{W\!F}}
		&\text{for any } f\in \mathcal{C}^2_b,
		\with
		\cL_{W\!F} f(x) :=  x\,(1 - x) 
		\lc \gamma_I\, \partial^2_{xx} f(x) - s\, \partial_{x}f(x)\rc,
		\\&
		N_t^f = 
		\LAg\nu_t\bv f\RAg 
		-\LAg\nu_0\bv f\RAg
		- \intO{t} \LAg\nu_u \bv \cL_{W\!F} f\RAg\, du
		+  \gamma_G\,
		\intO{t} \LAg\nu_u\bv 
		(\rho - \LAg\nu_u\bv \rho\RAg)\times  f
		\RAg\, du
\end{align*}
is a martingale with conditional quadratic variation~:
\begin{align*}
		&<N^f>_t = (\gamma_G)^2
		\intO{t} \Lbr\LAg\nu_u\bv f^2\RAg 
		- \LAg\nu_u\bv f\RAg^2 \Rbr\, du.
\end{align*}

There is actually an intermediate limit between these two, 
which will be the main focus of the current paper.
In this limit,
the fluctuations inside groups still play a role
while the fluctuations between groups are neglected
(rather in order to simplify the following analysis than for biological relevance) :
\begin{theo}
	\label{TS:MLim}
	Suppose
	that as $n, m\ifty$,
	we have the convergence of the rates 
	$\bar{\gamma}_I/n \rightarrow  \gamma_I$, 
	$\bar{\gamma}_I\, \bar{s} \rightarrow s$,
	while  
	$\limsup \bar{\gamma}_G < \infty$ 
	and $\{\bar{\gamma}_G\, \bar{r}(x)\}_{x\in [0,1]} \equiv \{r(x)\}_{x\in [0,1]}$ 
	is the same bounded measurable function for any $n, m$.
	Suppose that 
	$\mu^{m;n}_0$ is defined by assigning
	independently 
	the proportion of type C in each group
	according to the measure $\bar{\mu}^{m;n}_0$,
	where $\bar{\mu}^{m;n}_0 \rightarrow \mu_0$ as $m, n\ifty$.
	Then, 
	$\mu^{m;n}_t$ converges weakly in $D([0; T]; \M_1([0;1]))$
	to $\mu_t$, where
	$\mu_t$ is the unique solution to satisfy
	for any $f\in C^2_b$ the following equation	~:
	\begin{align*}
			&\partial_t\, \LAg\mu_t\bv f\RAg
			= \LAg \mu_t\bv \cL_{W\!F}f \RAg
			+  \LAg \mu_t\bv r\,f\RAg  
			- \LAg \mu_t\bv f\RAg \ltm \LAg \mu_t\bv r\RAg
			\mVg \qquad
			\mu_0 = \mu_0.
			\EQn{TS:eqCar}{}
	\end{align*}
\end{theo}

Since $\gamma_I$ is the only diffusion term left in this limit, 
we shall drop the subscript $I$ from now on.
The uniqueness of the solution is proved 
in the next Section 2.1
as part of Proposition \ref{TS:Char}, 
after we identify a more convenient way to describe it. 
The tightness of the above sequences 
and the fact that any limiting measure
is indeed a solution of \Req{TS:eqCar}
is deferred to Section \ref{TS:PMLim}.


\subsection{Definition of the solution of \Req{TS:eqCar} as a conditional law}
\label{TS:sec:MCNE}

Consider $X_t$ the $[0,1]$-valued solution of the following SDE, 
with initial condition $X_0\sim \mu_0$~:
\begin{align*}
		&dX_t := -s\, X_t\, (1-X_t) \, dt 
		+ \sqrt{2\gamma\, X_t\, (1-X_t)} \, dB_t.
		\EQn{TS:Xdef}{X}
\end{align*}
The existence and uniqueness of such a process can be found e.g. in in chapter 5.3.1 of \cite{D10}.

We will describe the solution of equation \Req{TS:eqCar} at time $t$
as the marginal distribution of $X_t$ with a Feynman-Kac penalization.
To relate the process 
to our results on convergence 
towards quasi-stationary distributions,
we will also represent this penalization 
as a conditioning upon survival 
of the stochastic process.

Since subtracting a constant to $r$ does not change the value of 
$\LAg \mu_t\bv r\,f\RAg  
- \LAg \mu_t\bv f\RAg \; \LAg \mu_t\bv r\RAg$,
and recalling that $r$ is bounded,
we can easily rewrite \Req{TS:eqCar}
in terms of $\rho(x) = \|r\|_\infty - r(x)$,
which is non-negative and bounded :
\begin{align*}
		&\partial_t\, \LAg\mu_t\bv f\RAg
		= \LAg \mu_t\bv \cL_{W\!F}f \RAg
		-  \LAg \mu_t\bv (\rho - \LAg \mu_t\bv \rho\RAg)\,f\RAg  
		\mVg \qquad
		\mu_0 = \mu_0
		\EQn{TS:eqRho}{}
\end{align*}
We then consider 
the following Feynman-Kac penalization~:
\begin{align*} 
		\textstyle  Z_t := \exp[- \intO{t} \rho(X_s)\,ds]
		\EQn{TS:Zexp}{Z}
\end{align*}

\begin{prop}
	\label{TS:Char}
	
	Define for each $t\ge 0$ 
	the probability measure $\mu_t$
	by :
	\begin{align*}
			\langle \mu_t\, \big| \, f\rangle 
			:= \mathbb{E}\left[ f(X_t)\, Z_t \right] 
			/\, \mathbb{E}\left[ Z_t \right],
			\frl{f\in \mathcal{C}([0,1])}.
			\EQn{TS:muT}{}
	\end{align*}
	$(\mu_t)_{t\ge 0}$ is the unique solution of equation $\Req{TS:eqCar}$
	(and equation \Req{TS:eqRho}).
\end{prop}

This penalization 
can then be interpreted as the probability
that the process has survived 
while confronted to a death rate of $\rho$,
conditionally on $(X_t)_{t\ge 0}$. 
More precisely, 
with $T_\partial$ an exponential r.v. with rate 1
that is independent from $X$,
we define the extinction time as~:
\begin{align*}
		&\ext := \inf\Lbr t\ge 0\pv 
		-\ln(Z_t) \ge T_\partial\Rbr,
		\EQn{TS:ext}{\ext}
\end{align*}

Clearly, $0$ and $1$ are absorbing for the dynamics of $X$.
We will also treat these fixation events 
as another kind of extinction.
The hitting times of $0$ and $1$ 
are denoted $\tau_0$ and $\tau_1$, and we consider any combination :
\begin{align*}
		\extO := \ext \wedge \tau_0, 
		\quad \extU := \ext \wedge \tau_1, 
		\quad \tau_{0,1} := \tau_0 \wedge \tau_1,
		\quad \extOU := \ext \wedge \tau_0 \wedge \tau_1.
		\EQn{TS:tauOU}{\extOU}
\end{align*}
The extinction rates of $\delta_0$, 
i.e.  $\rho_0 = \rho(0)$,
and $\delta_1$,
i.e. $\rho_1 = \rho(1)$
will play a crucial role 
in the long-time behavior of $\mu_t$.

\begin{prop}
	\label{TS:DecExt}
	With the above notations,
	we then define for any $t\ge 0$
	the probability measure $\mu_t$
	by :
	\begin{equation}\EQn{TS:decExt}{}
		\mu_t = x^0_t\, \delta_0 + x^1_t\, \delta_1
		+x^\xi_t\, \xi_t
	\end{equation}
	\begin{align*}
			&\text{ where }	 
			x^0_t := \dfrac{\mathbb{E}\left[ Z_{\tau_0}\,\exp[-\rho_0(t-\tau_0)]\, 
				; \, \tau_0 < t \right]}{\mathbb{E}\left[ Z_t \right]}, \qquad 
			x^1_t := \dfrac{\mathbb{E}\left[ Z_{\tau_1}\,\exp[-\rho_1(t-\tau_1)]\, 
				; \, \tau_1 < t \right]}{\mathbb{E}\left[ Z_t \right]},
			\\& \hcm{1}
			x^\xi_t := 
			\dfrac{ \mathbb{E}\left[ Z_t\, ; \, t < \tau_{0,1} \right]}{\mathbb{E}\left[ Z_t \right]}, \qquad
			\\& \hcm{1}
			\langle \xi_t\, \big| \, f\rangle 
			:= \dfrac{\mathbb{E}\left[ f(X_t)\, Z_t\, ; \, t < \tau_{0,1} \right] }
			{\mathbb{E}\left[ Z_t\, ; \, t < \tau_{0,1} \right]}
			= \mathbb{E}\left[ f(X_t)\, \vert \, t < \tau_{0,1,\partial} \right],
			\frl{f\in \mathcal{C}([0,1])}.
	\end{align*}
	$(\mu_t)_{t\ge 0}$ is the unique solution of equation \Req{TS:eqCar},
	which is equivalently given by :
	$\mu_t(dx) := \PR_{\mu_0}(X_t \in dx \bv t < \ext)$.
\end{prop}

%

\begin{rem}
	The solution of $\Req{TS:eqCar}$ 
	will thus generally be denoted $\mu_0 A_t$ 
	in the following statements.
	Expressing the dynamics in terms of an extinction rate 
	is done mainly to simplify notations with conditional laws. 
	We just adjusted the reference growth rate, 
	here $\|r\|_\infty$, to ensure that 
	the associated semi-group is sub-conservative.
\end{rem}
\begin{rem}
	In practice, 
	it means that one weights 
	specifically any potential trajectory 
	for the proportion of type C individuals inside a group.
	In order to obtain the dynamics that is typical while looking
	in the past of a uniformly sampled group,
	the weight of such trajectories
	is related to the mean number of lineages
	that are expected to follow this dynamics.
	For instance, 
	spending time where 
	the reproduction rate is high
	gives more opportunities 
	for at least one group to follow the trajectory until the end.
	In this view, 
	note that the solution to equation \Req{TS:Xdef} 
	is well-known to describe 
	the evolving proportion of an allele 
	under selection 
	in a population without any group
	(hence no selective effects between groups), 
	cf for instance \cite{Ew04}.
\end{rem}

\noindent
\textsl{Proof of Proposition \ref{TS:Char}:}
By the Ito formula, for any $f\in \C^2_b$~:
\begin{align*}
		&\E\lc f(X_t)\, Z_t \rc 
		= \LAg\mu_0\bv f\RAg
		+ \intO{t} \E\lc  \cL_{W\!F} f(X_s)\, Z_s \rc \, ds
		- \intO{t} \E\lc f(X_s)\,\rho(X_s)\, Z_s  \rc \, ds,
		\\&\hcm{1}
		\E\lc Z_t \rc 
		= 1 - \intO{t} \E\lc \rho(X_s)\, Z_s \rc \, ds,
\end{align*}
Thus~:
\begin{align*}
		&\partial_t\,\LAg \mu_t\bv f\RAg  
		= \dfrac{\E\lc \cL_{W\!F} f(X_t)\, Z_t \rc}{ \E\lc Z_t \rc  }
		-  \dfrac{\E\lc f(X_t)\,  \rho(X_t)\,Z_t \rc}{ \E\lc Z_t \rc  }
		+ \dfrac{\E\lc f(X_t)\, Z_t \rc}{ \E\lc Z_t \rc  } 
		\times \dfrac{\E\lc \rho(X_t)\, Z_t \rc}{ \E\lc Z_t \rc  } 
		\\&\hcm{2}
		= \LAg \mu_t\bv \cL_{W\!F} f \RAg
		+  \LAg \mu_t\bv r\,f\RAg  
		- \LAg \mu_t\bv f\RAg \times \LAg \mu_t\bv r\RAg.
\end{align*}
$(\mu_t)$ is indeed solution to equation \Req{TS:eqCar}.

Now, we turn to uniqueness.
Let $\bar{\mu}$ be a solution to equation \Req{TS:eqCar}, 
$P_t$ the  semi-group associated to $X_t$, 
the Wright-Fisher diffusion defined by \Req{TS:Xdef},
$f_0\in C^2_b([0,1])$, 
and for $0\le s\le t$~:
\begin{align*}
		&\hcm{1}
		\bar{n}_t : = \exp\lc \intO{t} \LAg \bar{\mu}_s\bv r \RAg\, ds \rc\mVg
		\hcm{1}
		f_s^t(x) = \bar{n_s}\times P_{t-s} f_0(x),
		\\& \text{so that~:}\hcm{2}
		\partial_s f_s^t(x) 
		:= \LAg \bar{\mu}_s\bv r \RAg\, f_s^t\  (x)  - \cL_{W\!F} f_s^t\  (x),
\end{align*}
\begin{align*}
		&\LAg \bar{\mu}_t\bv \bar{n_t}\, f_0 \RAg
		=\LAg \bar{\mu}_t\bv f_t^t \RAg
		:= \LAg \bar{\mu}_0\bv P_t\, f_0 \RAg
		+\intO{t} \big[
		\LAg \bar{\mu}_s\bv \cL_{W\!F}f_s^t \RAg
		+  \LAg \bar{\mu}_s\bv r\,f_s^t\RAg  
		- \LAg \bar{\mu}_s\bv f_s^t\RAg \times \LAg \bar{\mu}_s\bv r\RAg 
		\\&\hcm{5}
		+ \LAg \bar{\mu}_s\bv  \LAg \bar{\mu}_s\bv r \RAg \times f_s^t \RAg
		- \LAg \bar{\mu}_s\bv \cL_{W\!F} f_s^t \RAg  \big]
		\, ds,
		\\&\text{so that }
		\bar{\nu}_t(dx) := \bar{n}_t\, \bar{\mu}_t(dx) \quad
		\text{ solves }\quad
		\LAg \bar{\nu}_t\bv f_0 \RAg
		= \LAg \bar{\nu}_0\bv P_t\, f_0 \RAg
		+\intO{t}  \LAg \bar{\nu}_s\bv r\times P_{t-s}\,f_0\RAg \, ds.
		\EQn{TS:bnut}{\bar{\nu}}
\end{align*}

Recalling that we already have a solution $\mu_t$
defined through equation \Req{TS:muT}, 
we define similarly~: 
\begin{align*}
		n_t : = \exp\lc \intO{t} \LAg \mu_s\bv r \RAg\, ds \rc\mVg
		\hcm{.3}
		\nu_t(dx) := n_t\, \mu_t(dx).
		\EQn{TS:nuT}{\nu}
\end{align*}
As previously, $\nu$ is also solution to \Req{TS:bnut}, and we deduce :
\begin{align*}
		&|\LAg \nu_t- \bar{\nu}_t\bv f_0 \RAg|
		\le \intO{t}  |\LAg \nu_s- \bar{\nu}_s\bv r\times P_{t-s}\,f_0\RAg| \, ds
		\\&\hcm{2}
		\le 2\, \Ninf{f_0}\times \Ninf{r}\; \intO{t}   \NTV{\nu_s- \bar{\nu}_s} \, ds,
\end{align*}
with the convention :
$$\NTV{\nu}
= \sup_{f_0 \in C^2_b([0,1])}\dfrac{|\LAg \nu\bv f_0 \RAg|}{2\Ninf{f_0}}$$
Since this is true for any $f_0 \in C^2_b([0,1])$~: 
\begin{align*}
		& \NTV{\nu_t - \bar{\nu}_t}
		=\dfrac{|\LAg \nu_t- \bar{\nu}_t\bv f_0 \RAg|}{2\Ninf{f_0}}
		\le \Ninf{r}\; \intO{t}   \NTV{\nu_s- \bar{\nu}_s} \, ds.
\end{align*}
By Gronwall's Lemma (with the total variation uniformly bounded), 
this proves that $\nu_t = \bar{\nu}_t$ for any $t>0$.
Since $\mu_t$ (resp. $\bar{\mu}_t$) 
is deduced from $\nu_t$ (resp. $\bar{\nu}_t$) 
by renormalization, 
$\bar{\mu}_t = \mu_t$ for any $t>0$.
\epf
\\

\textsl{Proof of Proposition \ref{TS:DecExt}:}
First of all, we note that $\{t< \ext\} = \{\int_0^t \rho(X_s) ds < T_\partial\}$.
By the independence between $X$ and $T_\partial$,
$\PR(t < \ext \bv X) = Z_t$.
Thus :
\begin{align*}
		\mathbb{E}\left[ f(X_t)\pv t < \ext \right] 
		= \mathbb{E}\left[ f(X_t)\, Z_t \right] 
		\text{and in particular }	
		\PR(t < \ext)
		= \E(Z_t).
\end{align*}
This proves that for any $t\ge 0$
$\PR_{\mu_0}(X_t \in dx \bv t < \ext)$
defines a solution to equation \Req{TS:muT}.
By Proposition \ref{TS:Char}, it coincides with the solution with equation \Req{TS:eqCar}.

Since $0$ and $1$ are absorbing for the process $X$,
we have on the event  $\{\tau_0 < t\}$
(resp. $\{\tau_1 <~t\}$)
the fact that $Z_t := Z_{\tau_0}\, \exp[-\rho_0 (t- \tau_0)]$
and $X_t = 0$
(resp. $Z_t := Z_{\tau_0}\, \exp[-\rho_1 (t- \tau_1)]$ and $X_t =1$).
From this, Proposition \ref{TS:DecExt} is elementary.
\epf
\\

\subsection{Motivation for looking at $(\mu_t)$.$\quad$} 
\label{TS:sec:Motiv}

Our purpose in the following analysis of the long-time behavior of $\mu_t$
is to highlight common features
with the individual-based models $\mu^{m, n}_t$, 
for $m$ and $n$ reasonably large.
Note that contrary to $\mu_t$,
$\mu^{m, n}_t$ is a random process,
whose state space is the discrete grid :
\begin{equation}\EQn{TS:Mnm}{}
	\M^{m, n}([0,1]) 
	:= \{\mu \in \MOne 
	\bv \mu(\medcup_{0\le k\le n}\{k/n\}) = 1\pv
	\frl{0\le k\le n} m \ltm \mu(\{k/n\}) \in [\![0, m]\!]
	\}.
\end{equation}
So we generally expect a greater diversity of possible scenarii
for the  individual-based models.
The mathematical description of its limiting behavior
is however not very informative unless one gets precise quantitative estimates,
as we shall see in the next lemma.
\begin{lem}
	\label{TS:QSD.MN}
	For any $m, n \ge 1$,
	$\delta_0$ and $\delta_1$
	are the only absorbing points of $\mu^{m, n}_t$.
	Denoting $\tau_{0,1}^{m, n}$ this absorption time,
	there exists a unique associated QSD $A^{m, n} \in \M_1(\M^{m, n}([0,1]))$
	and a unique capacity of survival $H^{m, n} \in L_\infty(\M^{m, n}([0,1]))$
	with extinction rate $\rho_A$.
	Moreover, there exists $C, \chi$ such that 
	the following convergence results hold :
	\begin{align*}
			&\frlq{\mu_0, \mu_1 \in \M^{m, n}([0,1]) \setminus \{\delta_0, \delta_1\}}
			\\
			&0\le \PR_{\mu_0}(\Ex{s\ge 0} \mu^{m, n}_s = \delta_0)
			- \PR_{\mu_0}(\mu^{m, n}_t = \delta_0)
			\le C \exp(- \rho_A\, t),
			\\& 0\le \PR_{\mu_0}(\Ex{s\ge 0} \mu^{m, n}_s = \delta_1)
			- \PR_{\mu_0}(\mu^{m, n}_t = \delta_1)
			\le C \exp(- \rho_A\, t),
			\\
			&|\exp(\rho_A\, t) \PR_{\mu_0}(\mu^{m, n}_t = \mu_1)
			- H^{m, n}(\mu_0)\ltm A^{m, n} (\{\mu_1\})| \le C \exp(- \chi\, t).
	\end{align*}
\end{lem}
The proof of this lemma is deferred to the Appendix.
In this proof, the parameters $C, \rho_A$ and $\chi$
strongly depend on $m$ and $n$. 
So there is no clear dependency on the parameters $r$ and $s$ in this lemma,
unless one gets by other means precise estimates of $C, \rho_A$ and $\chi$
as well as $\PR_{\mu_0}(\Ex{s\ge 0} \mu^{m, n}_s = \delta_0)$ 
and $\PR_{\mu_0}(\Ex{s\ge 0} \mu^{m, n}_s = \delta_1)$. 
The study of $A^{m, n}$ and $H^{m, n}$ is also much more complicated 
than the one of $\alpha$ and $\heig$ of Proposition \ref{TS:p:QSD01}
because of the state space of the former, namely $\M^{m, n}([0,1])$,
is much more complicated than the one of the latter, namely $[0,1]$.

The description of $\mu_t$  
as the conditional distribution $\mu_0 A_t$ 
provides us with very efficient tools 
for asymptotic results,
as we shall see in Section \ref{TS:sec:PfQSD}.
First, this sheds light on specific conditions 
separating different modes of convergence.
The limit is also identified uniquely 
for any set of parameters and initial condition.
We finally obtain exponential convergence in total variation,
whose rate can be somewhat identified and interpreted.

This is why we study $(\mu_t)$
in order to describe the dynamics of $(\mu^{m, n}_t)$,
as a complement to the study provided in \cite{LM15}
of the solution $(\pi_t)$ of equation \Req{TS:piT}.
A discussion is presented in Section \ref{TS:sec:Disc}
to see how much information we can actually infer 
from the long term behavior of $(\mu_t)$
regarding the one of $(\mu^{m, n}_t)$
and describe the main limitations that we face.

The solution $(\nu_t)$ of the martingale problem
given by equation \Req{TS:nuT} 
is likely to behave more closely to $(\mu^{m, n}_t)$.
Yet, we expect a poorly informative limiting behavior
similar to the one of Lemma \ref{TS:QSD.MN}
with a proof that seems too technical for now.
Furthermore,
regarding the limitations that we plan to describe in Section \ref{TS:sec:Disc},
we shall have similar issues in the connection of 
$(\mu^{m, n}_t)$ to $\nu_t$.

To gain some more perspective on our analysis
on the long time behavior 
of the solution $\mu_t$ of equation \Req{TS:eqCar},
we have used a numerical approximation.
What we present in Section \ref{TS:sec:Disc} 
relies on the intuitions provided by these simulations.
The subject would 
require a much more complete study
to be more quantitative,
but these simulations already provide 
illustrations that our convergence results
can be informative.
The simulations are also helpful 
to highlight some of the above-mentioned limitations
in relating the limiting behavior of $(\mu^{m, n}_t)$
and of $(\mu_t)$.

\section{QSDs and exponential convergence}
\label{TS:sec:QSD}
\setcounter{eq}{0}

We know from Proposition \ref{TS:DecExt}
that the quasi-stationary distributions for $X$ with extinction at time $\ext$
correspond exactly to the initial conditions 
for which the solution of equation \Req{TS:eqRho} is constant in time.
Such a QSD will be called stable 
if it is the quasi-limiting distribution
for any initial condition 
that is close enough in total variation distance.

Since $0$ and $1$ are absorbing states for $X$,
$\delta_0$ and $\delta_1$ necessarily belong to these QSDs,
with extinction rate respectively $\rho_0$ and $\rho_1$.

We define the semi-groups associated to our different extinctions ~: 
\begin{align*}
		&
		\mu P_t (dx) 
		:= \PR_\mu (X_t \in dx \pv t <\ext)\mVg
		\hcm{1.5}
		\mu A_t (dx) 
		:= \PR_\mu (X_t \in dx \bv t <\ext)
		\\&
		\mu P^{01}_t (dx) 
		:= \PR_\mu (X_t \in dx \pv t <\extOU)\mVg
		\hcm{.8}
		\mu A^{01}_t (dx) 
		:= \PR_\mu (X_t \in dx \bv t <\extOU)
		\\&
		\mu P^{1}_t (dx) 
		:= \PR_\mu (X_t \in dx \pv t <\extU)\mVg
		\hcm{1.1}
		\mu A^{1}_t (dx) 
		:= \PR_\mu (X_t \in dx \bv t <\extU)
\end{align*} 

\begin{prop}
	\label{TS:p:QSD01}
	There exists a unique QSD $\alpha\in \MoneM$ 
	and a unique capacity of survival $\heig$
	associated to the extinction $\extOU$.
	With the associated extinction rate $\rho_\alpha$, it means first that~:
	\begin{align*}
			\frlq{t>0}
			&\alpha P^{01}_t (dx) 
			= \exp[-\ra\, t]\, \alpha(dx)
			\mVg 
			P^{01}_t \heig = \exp[-\ra\, t]\, \heig
	\end{align*}
	Moreover, 
	for any $\mu \in \MoneM$ and $x\in (0,1)$~:
	\begin{align*}
			&\alpha(dx) = \lim_{t\ifty} \mu A^{01}(dx)\mVg \quad
			\heig(x) = \lim_{t\ifty} h_t(x)
			\\&
			\where h_t(x) := \exp[\ra\, t]\,\PR_x(t<\extOU).
			\EQn{TS:ht}{h_t}
	\end{align*}
	The convergence are uniformly exponential,
	in the sense that there exists $\chi, C >0$ such that~:
	\begin{align*}
			&\frlq{\mu \in \MoneM}
			\NTV{\mu A^{01}_t
				- \alpha}
			\vee \Ninf{h_t - h}
			\le C \, \exp[-\chi\, t].
			\EQn{TS:CVal}{\alpha}
			%
			\\& 
			\text{ in particular }\quad
			\Ninf{\bar{\heig}}
			:= \sup\!_{\{x\in (0,1),\, t>0\}}
			\exp[\rho_\alpha\, t]\, \PR_x(t<\extOU) < \infty
			\EQn{TS:etaB}{\Ninf{\bar{\heig}}}
	\end{align*}
	Moreover, for any $n\ge 2$,
	$h$ is lower-bounded by a positive constant on $[1/n, 1- 1/n]$.
\end{prop}
We show in the following Subsections that the long-time behavior
of the process with only the local extinction rate
depends mainly on $\rho_\alpha,$ $\rho_0$ and $\rho_1$.
In the convergences that follow, 
we will often have uniform bounds for probability measures 
belonging for some $n \ge 1$ and $\xi \in (0,1)$ to~:
\begin{align*}
		\MnxTS 
		:= \Lbr \mu\in \MOne\bv
		\mu[1/n, 1] \ge \xi\Rbr\mVg
		\medcup_{n, \xi} \MnxTS  
		= \MOne \setminus \{\delta_0\}.
		\EQn{TS:Mnx}{\MnxTS}
\end{align*}
\begin{align*}
		&\text{or in}\quad
		\MnxTS^{0,1} 
		:= \Lbr \mu\in \MOne\bv
		\mu[1/n,\, 1-1/n] \ge \xi\Rbr\mVg
		\quad(n\ge 3, \xi >0)
		\EQn{TS:MnxOU}{\MnxTS^{0,1}}
		\\&
		\medcup_{n, \xi} \MnxTS^{0,1}  
		= \MOne \setminus 
		\{x\, \delta_0 + (1-x)\, \delta_1\bv x \in [0,1]\}.
\end{align*}

But first, the following Lemma provides some elementary properties
of the extinction rate in terms of the function $\rho$.
\begin{lem}
	\label{TS:ExtRho}
	Assume that $\mu_\infty$ is a QSD of the Markov process $X$
	for the extinction time $\ext$ defined by \Req{TS:ext}.
	Then, its extinction rate is given by $\LAg \mu_\infty \bv \rho\RAg$.
	
	The QSD $\alpha$  given in Proposition \ref{TS:p:QSD01}
	satisfies :
	$$\LAg \alpha\bv \rho\RAg = \rho_\alpha\, \PR_\alpha(\extOU = \ext) < \rho_\alpha.$$
\end{lem}

\paragraph{Proof of Lemma \ref{TS:ExtRho}}
Let $\lambda$ be the extinction rate of the QSD $\mu_\infty$. 
We thus know that for any $t>0$ :
\begin{align*}
		\lambda 
		&= \frac{-1}{t} \log \PR_{\mu_\infty}(t< \ext)
		=\frac{-1}{t} \log \E_{\mu_\infty}(Z_t)
		\\&=\frac{-1}{t} \log \E_{\mu_\infty}(\exp[-\int_0^t \rho(X_s) ds])
\end{align*}
We can simply look at the limit of this expression as $t$ tends to 0
to deduce Lemma \ref{TS:ExtRho}. 
With the Fubini Theorem,
the expression can also be identified with fixed $t$, because :
\begin{align*}
		&\E_{\mu_\infty}(\exp[-\int_0^t \rho(X_s) ds])
		= 1 - \int_0^t ds\, \E_{\mu_\infty}( \rho(X_s)\,\exp[-\int_0^s \rho(X_u) du])
		\\&\hcm{1}
		= 1 - \int_0^t ds\, \LAg \mu_\infty P_t \bv \rho\RAg
		= 1 - \LAg \mu_\infty \bv \rho\RAg \ltm \dfrac{1-\exp{-\lambda t}}{\lambda}.
\end{align*}
Concerning $\alpha$,
we can exploit Theorem 2.6 in \cite{coll}, 
that proves that the exit state is independent from the exit time 
when the initial condition is a QSD, 
with an exponential law for the exit time.
This implies notably that for any $t\ge 0$~:
\begin{equation}\EQn{TS:ExtInd}{}
	\E_\alpha(\exp[\rho_\alpha \extOU]\pv \extOU = \ext \le t)
	= \rho_\alpha^t\, \PR_\alpha(\extOU = \ext).
\end{equation}
Then, we can prove :
\begin{align*}
		\LAg \alpha\bv \rho\RAg
		&\textstyle
		= \lim_{t\rightarrow 0} 
		\frac{1}{t} \PR_\alpha(T_\partial < \int_0^t \rho(X_s) ds)
		= \lim_{t\rightarrow 0} 
		\frac{1}{t} \PR_\alpha(\ext \le t < \tau_{0,1})
		\\&\textstyle
		= \lim_{t\rightarrow 0} 
		\frac{1}{t} \E_\alpha(\exp[\rho_\alpha \extOU]\pv \extOU = \ext \le t).
\end{align*}
By \Req{TS:ExtInd},
this proves $\LAg \alpha\bv \rho\RAg = \rho_\alpha\, \PR_\alpha(\extOU = \ext) < \rho_\alpha$.

%
%
%
%

\subsection{$\rho_1 < \rho_0 < \rho_\alpha$ : 
	Group selection favoring one allele with a quick fixation}


\begin{prop}
	\label{TS:p:r1r0ra}
	Assume that $\rho_1 < \rho_0 < \rho_\alpha$.
	$\delta_1$ is then the only stable QSD, 
	with convergence rate $\rho_0 -\rho_1$,
	i.e. :
	\begin{align*}
			\frl{n\ge 1}\frl{\xi >0}
			\Exq{C_{n, \xi}>0}
			\frlq{\mu \in \MnxTS}
			\NTV{\mu A_t - \delta_1}
			\le C_{n, \xi} \, \exp[-(\rho_0 -\rho_1)\, t].
	\end{align*}
\end{prop}

In order to obtain this convergence,
we considered the dynamics of $\mu A_t$ restricted on $[0, 1)$,
which happens to be given by $\mu A^1_t$ 
because the process is fixed at $1$ until its extinction,
once $1$ is reached.
Looking at the dynamics of $\mu A^1_t$,
we deduce an additional level of convergence
for the disappearance of polymorphic groups :
\begin{prop}
	\label{TS:p2:r1r0ra}
	Assume that $\rho_0 < \ra$. Then, there exists $C>0$ such that~:
	\begin{align*}
			\frlq{\mu \in \MOne\setminus\{\delta_1\}}
			\NTV{\mu A^{1}_t
				- \delta_0}
			\le C \, \exp[-(\rho_\alpha -\rho_0)\, t].
	\end{align*}
\end{prop}
With the notations of Proposition \ref{TS:DecExt},
the two previous propositions imply that :
\begin{align*}
		x_t^0\le 1 - x_t^1\le C_{n, \xi} \, \exp^{-(\rho_0 -\rho_1)\, t}
		\pv \quad 
		x_t^\xi \le C \,x_t^0\, \exp^{-(\rho_\alpha -\rho_0)\, t}
		\le C\ltm  C_{n, \xi}\, \exp^{-(\rho_\alpha -\rho_1)\, t}.
\end{align*}

\begin{rem}
	In this limiting model,
	whatever the selection effects inside the groups,
	the selective effects at group level
	favoring any of the pure groups 
	always dominate in the long run.
	The convergence to the pure C groups population
	happens in total variation,
	with at the end an exponential rate of convergence.
	This rate is given
	by the competition between pure groups.
\end{rem}

\begin{rem} 
	\label{TS:relQSD}
	We shall precise in Section \ref{TS:sec:Disc}
	and more specifically in Section \ref{TS:sec:pureSel}
	the limits of this description,
	in particular when one wishes to relate it
	in terms of the individual-based model.
	Even if pure C groups happen to dominate 
	in the long run,
	expecting an exponential convergence rate
	might be misleading~:
	the initial proportion of pure $C$ groups 
	may be so small 
	that pure $D$ groups would be dominant 
	for a very long time.
	Some illustrations obtained by simulations
	of such a case
	are  given in Figures \ref{TS:FD26}-\ref{TS:QSD26}
	around Subsection \ref{TS:sec:pureSel}.
	The main quantities of interest are then 
	the time needed
	for the competition between groups 
	to compensate the initial domination by pure $D$ groups,
	assuming that it happens,
	and the probability that this transition actually occurs 
	for the IBM.
\end{rem}

\subsection{$\rho_1 < \rho_\alpha < \rho_0$ : 
	Group selection favoring one allele with a slow fixation}

\begin{prop}
	\label{TS:p:r1rar0}
	Assume that $\rho_1 < \rho_\alpha < \rho_0$.
	Then, $\delta_1$ is again the only stable QSD, 
	with convergence rate $\rho_\alpha -\rho_1$,
	i.e. :
	\begin{align*}
			\frl{n\ge 1}\frl{\xi >0}
			\Exq{C_{n, \xi}>0}
			\frlq{\mu \in \MnxTS}
			\NTV{\mu A_t - \delta_1}
			\le C_{n, \xi} \, \exp[-(\rho_\alpha -\rho_1)\, t].
	\end{align*}
\end{prop}

Again, we have an additional level of convergence,
and the quasi-equilibrium is precisely described in terms
of the polymorphic quasi-stationary distribution :
\begin{prop}
	\label{TS:p2:r1rar0}
	Assume that $\ra < \rho_0$. Then :
	\begin{align*}
			&\Exq{\chi^1 > 0}
			\frl{n\ge 1}\frl{\xi >0}
			\Exq{C_{n, \xi}>0}
			\frlq{\mu \in \MnxTS\setminus \{\delta_1\}}
			\\&\hcm{2}
			\NTV{\mu A^{1}_t
				- \alpha_1}
			\le C_{n, \xi} \, \exp[-\chi^1\, t],
			\EQn{TS:alU}{}
	\end{align*}
	where the QSD $\alpha_1$ has extinction rate $\rho_\alpha$
	and is given as 
	$\alpha_1 = y_0\, \delta_0 + y_\alpha\, \alpha$
	with the relations~:
	\begin{align*}
			&\hcm{2} \dfrac{y_0}{y_\alpha} 
			= \dfrac{\ra\times \PR_\alpha(\tau_0 = \extOU) }{
				(\rho_0 -\rho_\alpha)}\mVg
			\qquad
			y_0 + y_\alpha = 1
			\EQn{TS:y0}{}
			\\&
			\text{ and thus }
			y_\alpha 
			:= \dfrac{  (\rho_0 -\rho_\alpha)}
			{\rho_0 -\ra\times \PR_\alpha(\extU =\extOU) }\mVg
			\quad
			y_0
			:= \dfrac{\ra\times \PR_\alpha(\tau_0 = \extOU) }
			{\rho_0 -\ra\times \PR_\alpha(\extU =\extOU) }.
	\end{align*}
	Moreover, 
	we know the associated capacity of survival $\heig^1 := \heig / y_\alpha$
	($\heig^1(0) = 0$) and :
	\begin{align*}
			&\frl{n, \xi}
			\Exq{C_{n, \xi}>0}
			\frlq{\mu \in \MnxTS}
			\\& \hcm{2}
			|\exp[\ra\, t]\, \PR_\mu(t<\extU)
			- \LAg \mu\bv \heig^1\RAg | 
			\le C_{n, \xi} \, \exp[-\chi^1\, t]
			\EQn{TS:CVeta1}{}
			\\& \hcm{1}
			\text{ and }\quad
			\Ninf{\bar{\heig}^1} 
			:= \sup\!_{\{x\in [0,1),\, t>0\}}
			\exp[\rho_\alpha\, t]\, \PR_x(t<\extU) < \infty.
			\EQn{TS:etaB1}{\Ninf{\bar{\heig}^1} }
	\end{align*}
\end{prop}

With the notations of Proposition \ref{TS:DecExt},
the two previous propositions notably imply that :
\begin{align*}
		x_t^0 + x_t^\xi = 1 - x_t^1\le C_{n, \xi} \, \exp^{-(\rho_0 -\rho_1)\, t}
		\pv \quad 
		\frac{x_t^\xi}{x_t^0}  \rightarrow \frac{y_\alpha}{y_0}.
\end{align*}

\begin{rem}
	Some illustrations obtained by simulations
	of this situation
	are  given in Figures \ref{TS:FD21}-\ref{TS:FDL21}
	around Subsection \ref{TS:InterG}.
	The maintenance of pure $D$ groups 
	in the population is here mainly due 
	to the fixation of polymorphic groups,
	so that their proportion relative to polymorphic groups 
	tends to stabilize, while both vanish.
	It does not affect the asymptotic profile 
	of these polymorphic groups in this limit.
	Only the proportion of these groups is adjusted :
	extinction of pure $D$ groups
	exactly compensates their generation 
	by fixation of polymorphic ones.
\end{rem}

\begin{rem}
	If we were to include 
	a small effect of neutral replacements of groups, 
	or consider the individual-based models,
	this size reduction
	would imply a higher rate of these fluctuations 
	for the polymorphic profile,
	as compared to having fixation 
	implying extinction of the groups.
\end{rem}

\begin{rem}
	It may happen that 
	the polymorphic QSD actually emerges 
	after a long domination of the other pure groups $D$.
	First results of simulations indicate 
	that this quasi-equilibrium 
	may not even be noticeable
	when looking at the population as a whole
	since it can emerge at almost the same time 
	as the pure groups,
	cf Section \ref{TS:sec:pureSel}.
	Some illustrations obtained by simulations
	of such a situation
	are also given in Figures \ref{TS:FDL20}, \ref{TS:FDLC20}, \ref{TS:FD21}, \ref{TS:FDL21}
	in Section \ref{TS:sec:Disc}.
\end{rem}

\subsection{$\rho_1 < \rho_0 = \rho_\alpha$ : 
	Group selection favoring one allele with a critical fixation rate}

\begin{prop}
	\label{TS:p:r1r0a}
	Assume that $\rho_1 < \rho_0 = \rho_\alpha$.
	Then, $\delta_1$ is again the only stable QSD, 
	with convergence rate $\rho_0 -\rho_1$,
	i.e. :
	\begin{align*}
			\frl{n\ge 1}\frl{\xi >0}
			\Exq{C_{n, \xi}>0}
			\frlq{\mu \in \MnxTS}
			\NTV{\mu A_t - \delta_1}
			\le C_{n, \xi}\ltm (1+t) \, \exp[-(\rho_0 -\rho_1)\, t].
	\end{align*}
\end{prop}

For the next level of convergence,
$0$ is still dominant, 
yet the proportion of polymorphic states 
is vanishing in comparison only at rate $1/t$ :
\begin{prop}
	\label{TS:p2:r1r0a}
	Assume that $\rho_0 = \rho_\alpha$.
	Then, $\delta_0$ is the only QSD with extinction $\extU$. 
	Moreover~:
	\begin{align*}
			&\Exq{t_\vee, C>0}
			\frl{t\ge t_\vee}
			\frlq{\mu \in \MOne}
			\hcm{0.8}
			\NTV{\mu A^{1}_t
				- \delta_0}
			\le C  / (1+t),
			\\&
			\frl{n, \xi} 
			\Exq{t_{n,\xi}, C_{n,\xi} >0}
			\frl{t\ge t_{n,\xi}}
			\frlq{\mu \in \MnxTS^{0,1}}
			\NTV{\mu A^{1}_t
				- \delta_0} 
			\ge C_{n,\xi} / t.
	\end{align*}
\end{prop}

With the notations of Proposition \ref{TS:DecExt},
the two previous propositions notably imply that provided $x_0^\xi >0$~:
\begin{align*}
		x_t^0 \le  1 - x_t^1\le C_{n, \xi} (1+t)\, \exp^{-(\rho_0 -\rho_1)\, t}
		\pv \quad 
		\frac{x_t^\xi}{x_t^0} = O_{t\ifty}(1/t).
\end{align*}

\begin{rem}
	In this critical case,
	pure $D$ groups happen to dominate 
	the transient dynamics,
	but only because of an asymptotically linear increase
	of their proportion.
	This increase is actually due to the fixation 
	of polymorphic groups :
	these groups, 
	still asymptotically distributed as $\alpha$,
	act as a generator of pure $D$ groups
	that stay conserved.
	Asymptotically, 
	the selective pressure is indeed uniform 
	between pure $D$ groups
	and the polymorphic groups as a whole.
\end{rem}

\begin{rem}
	The same issue of relevance
	as in Remark \ref{TS:relQSD} may be noted :
	the polymorphic QSD 
	might only be reached once pure $C$ groups 
	have already appeared non-negligible.
\end{rem}

\subsection{$\rho_0 = \rho_1 < \rho_\alpha$ : 
	selective effects at group level favoring fixation}

Again, the convergence rate of the distribution 
is given by the competition between the pure groups 
and the polymorpic QSD.
There is in this case
a strong dependency on the initial condition 
regarding the final equilibrium.

\begin{prop}
	\label{TS:p:r10ra}
	Assume that $\rho_0 = \rho_1 < \rho_\alpha$.
	Then, any convex combination of  $\delta_0$ and $\delta_1$ 
	is a QSD, with extinction rate  $\rho_1$.
	The convergence still happens 
	with convergence rate $\rho_\alpha -\rho_1$,
	i.e. :
	\begin{align*}
			\Exq{C>0}
			\frl{\mu \in \MOne}
			\Exq{x \in [0,1]}
			\NTV{\mu A_t - (x\,\delta_0 + (1-x)\,\delta_1)}
			\le C \, \exp[-(\rho_\alpha -\rho_1)\, t].
	\end{align*}
	Moreover, the proportion $x$ for the limiting QSD is :
	\begin{align*}
			x(\mu)
			:= \E_\mu\lc \exp(\rho_1\, \extOU)\pv 
			\extOU = \tau_0\rc
			/ \E_\mu\lc \exp(\rho_1\, \extOU)\pv 
			\extOU = \tau_{0,1}\rc.
			\EQn{TS:xmu}{x(\mu)}
	\end{align*}
\end{prop}

With the notations of Proposition \ref{TS:DecExt},
this implies that for $x$ depending on the initial condition :
\begin{align*}
		|x_t^0 - x| \le C\, \exp^{-(\rho_\alpha -\rho_1)\, t}
		\pv \quad 
		|x_t^1 - 1 + x| \le C\, \exp^{-(\rho_\alpha -\rho_1)\, t}
		\pv \quad 
		|x_t^\xi| \le C\, \exp^{-(\rho_\alpha -\rho_1)\, t}
\end{align*}
The next level of convergence (with extinction $\extOU$)
is the already known 
convergence to $\alpha$ at exponential rate.
\\

\begin{rem}
	Here, the selective effects at group level are allowed 
	to dictate the dynamics of polymorphic groups.
	For a polymorphic (intermediate) initial condition,
	it may well happen
	(depending on $r$ and $\gamma$)
	that $x(\mu)$ is close to 0.
\end{rem}

\subsection{$\rho_\alpha < \rho_0\wedge \rho_1$ : 
	selective effects at group level strongly favoring polymorphism}


\begin{prop}
	\label{TS:p:rar10}
	Assume that $\rho_\alpha < \rho_0\wedge \rho_1 := \rho$.
	Then, there is only one stable QSD $\alpha^{0,1}$, 
	with convergence rate $\rho -\rho_\alpha$,
	i.e. :
	\begin{align*}
			\frl{n\ge 1}\frl{\xi >0}
			\Exq{C_{n, \xi}>0}
			\frlq{\mu \in \MnxTS^{0,1}}
			\NTV{\mu A_t - \alpha^{0,1}}
			\le C_{n, \xi} \, \exp[-(\rho_\alpha -\rho)\, t],
	\end{align*}
	where $\alpha^{0,1}$ has extinction rate $\rho_\alpha$
	and is given as 
	$\alpha^{0,1} = y_0\, \delta_0 +  y_1\,\delta_1  + y_\alpha\, \alpha$
	with~:
	\begin{align*}
			&
			\dfrac{y_0}{y_\alpha} 
			= \dfrac{\ra\times \PR_\alpha(\tau_0 = \extOU) }{
				(\rho_0 -\rho_\alpha)}\mVg
			\quad
			\dfrac{y_1}{y_\alpha} 
			= \dfrac{\ra\times \PR_\alpha(\tau_1 = \extOU) }{
				(\rho_1 -\rho_\alpha)}\mVg
			\EQn{TS:a01}{\alpha_{01}}
	\end{align*}
	and of course $y_0 + y_1 + y_\alpha = 1$.
	Regarding the capacity of survival $\heig^{0,1}$ :
	\begin{align*}
			\color{blue}
			&\frl{n, \xi}
			\Exq{C_{n, \xi}>0}
			\frlq{\mu \in \MnxTS}
			\\& \hcm{2}
			|\exp[\ra\, t]\, \PR_\mu(t<\extU)
			- \LAg \mu\bv \heig^{0,1}\RAg | 
			\le C_{n, \xi} \, \exp[-\chi^1\, t]
			\EQn{TS:CVeta01}{}
			\\& \hcm{1}
			\text{ and }\quad
			\Ninf{\bar{\heig}^{0,1}} 
			:= \sup\!_{\{x\in [0,1],\, t>0\}}
			\exp[\rho_\alpha\, t]\, \PR_x(t<\ext) < \infty.
			\EQn{TS:etaB01}{\Ninf{\bar{\heig}^1} }
	\end{align*}
	where for $x \in (0,1)$, $\heig^{0,1}(x) = h(x)/y_\alpha$
	while $\heig^{0,1}(0) = \heig^{0,1}(1) =0$.
\end{prop}

If $\rho_1 < \rho_0$, for any initial condition
$\mu = x\, \delta_0 + (1-x)\, \delta_1$ with $x\in(0,1)$,
$\mu A_t$ converges at rate $\rho_0 - \rho_1$ to $\delta_1$
as $t\ifty$.

If $\rho_1 = \rho_0$, then any such distribution is a QSD
with the extinction rate $\rho_0$.
\\

\begin{rem}
	Like in Proposition \ref{TS:p2:r1rar0},
	the QSD $\alpha^{0,1}$ is actually obtained
	by the stabilization of the profile
	of polymorphic groups (towards $\alpha$),
	then by the compensation 
	between the fixation of pure groups by $\alpha$
	and their extinction.
	Since $\heig^{0,1}$ is null at $0$ and $1$ 
	(it vanishes also at their vicinity),
	we see that 
	the contribution of fixed groups
	to the survival of the population becomes negligible.
	These pure groups are in fact driven 
	by the polymorphic groups, 
	as one could see from their lineages :
	it would not take long to come back 
	to polymorphic ancestors.
\end{rem}

\begin{rem}
	Considering an individual selection 
	depending on the frequency in the group, 
	we could easily extend our model 
	to describe the case of a balancing selection.
	In such an extension,
	as soon as fixation is not too exceptional,
	selective effects at group level are however needed 
	to maintain polymorphism
	without transmission
	between groups.
\end{rem}

\subsection{$\rho_1 = \rho_\alpha < \rho_0$ : 
	critical vanishing of the polymorphic QSD}

\begin{prop}
	\label{TS:p:r1ar0}
	Assume that $\rho_1 = \rho_\alpha < \rho_0$.
	Then, $\delta_1$ is again the only stable QSD, 
	yet the convergence is not exponential,
	and more precisely~:
	\begin{align*}
			&\hcm{1.2}
			\Exq{C>0}
			\frlq{\mu \in \MOne}
			\NTV{\mu A_t
				- \delta_1}
			\le C  / (1+t),
			\\& \frl{n\ge 2}
			\frlq{\xi>0}
			\Ex{c_{n, \xi}>0}
			\frlq{\mu \in \MnxTS^{0,1}}
			\\&\hcm{3}
			\NTV{\mu A_t
				- \delta_1}
			\ge c_{n, \xi}  / (1+t).
	\end{align*}
\end{prop}

For the next level of convergence, 
we refer to Proposition \ref{TS:p2:r1rar0}.
\\

\begin{rem}
	This case corresponds 
	to a very specific compensation 
	of the parameters,
	where selective effects at group level exactly compensate
	for the fixation events.
	The rate of convergence is slow,
	because it is driven by the polymorphic groups becoming negligible
	as compared to the fixed lineages that they generate.
\end{rem}

\subsection{$\rho_0 = \rho_1 = \rho_\alpha$ : 
}

This case is the most counter-intuitive, 
since any polymorphic component 
in the initial distribution 
imposes a predictable final equilibrium 
without polymorphism.

\begin{prop}
	\label{TS:p:r10a}
	Assume that $\rho_0 = \rho_1 = \rho_\alpha$.
	Then, any convex combination of  $\delta_0$ and $\delta_1$ 
	is a QSD, with extinction rate  $\rho_1$.
	They are the only ones, 
	and among them, only one is stable~:
	\begin{align*}
			\frlq{n, \xi}
			\Exq{C_{n,\xi}>0}
			\frlq{\mu \in \MnxTS^{0,1}}
			\NTV{\mu A_t - (x\,\delta_0 + (1-x)\,\delta_1)}
			\le C_{n,\xi} / (1+t).
	\end{align*}
	where the proportion $x$ for the limiting QSD is :
	\begin{align*}
			x
			:= \PR_\alpha(\tau_0 = \extOU) / \PR_\alpha(\tau_{0,1} = \extOU).
	\end{align*}
\end{prop}

\begin{rem}
	The distribution inside the interval vanishes so slowly
	that its flux to $0$ and $1$ governs the final equilibrium
	(with a much quicker stabilization of $\mu A^{0, 1}_t$ to $\alpha$).
\end{rem} 

\subsection{Limits of the parameters}

\begin{prop}
	\label{TS:prop:SI}
	Given any $s>0$ and any bounded function $r$, 
	$\lim_{\gamma \rightarrow \infty} \rho_\alpha(\gamma) = + \infty$.
\end{prop}

\begin{prop}
	\label{TS:prop:rI}
	Given any $\gamma >0$, $s\ge 0$, 
	and a continuous (and negative) function $r^0$ with its maximum 
	only in the interior of $(0,1)$,
	there exists a critical value $R_\vee >0$ such that
	for any $R>R_\vee$ and considering the system with $r = R\, r^0$,
	it holds $\rho_\alpha < \rho_0\wedge \rho_1$.
\end{prop}
Polymorphism is maintained 
by any sufficiently large selective effects at group level favoring it.

\begin{prop}
	\label{TS:prop:r0}
	Conversely,
	given any $\gamma >0$, $s\ge 0$, 
	and a bounded function $r^0$,
	there exists a critical value $R_\wedge >0$ such that
	for any $R<R_\wedge$ and considering the system with $r = R\, r^0$,
	it holds $\rho_0\wedge \rho_1< \rho_\alpha$.
\end{prop}
When the selective effects at group level is too small, 
polymorphism cannot maintain itself.
\\

\begin{rem}
	One could expect $\rho_\alpha(\gamma)$ 
	to be first a decreasing function of $\gamma$
	and then increasing. 
	Yet, it seems not to hold true for any general $r$.
	Think for instance of two types of equilibria that compete inside $(0,1)$,
	i.e. $r$ with two localized modes,
	with a specific optimal value $\gamma_1 < \gamma_2$ for each.
	In the range $\gamma = \gamma_1$ to $\gamma = \gamma_2$,
	the QSD shifts from the first mode, 
	where the extinction is becoming much larger as $\gamma$ increases,
	to the second mode, 
	where such increase is much less significant.
	It can happen if there is a very strong mode of $r$ close to a border,
	that is responsible for the first equilibrium.
	We may thus observe  $\rho_\alpha(\gamma) > \rho_\alpha(\gamma_1)\vee \rho_\alpha(\gamma_2)$
	for $\gamma \in (\gamma_1, \gamma_2)$,
	which contradicts the predicted profile of $\rho_\alpha$.
\end{rem}

We conjecture that
$\lim_{\gamma \rightarrow 0} \rho_\alpha(\gamma) 
= \infty$ also holds for any $s>0$ and any bounded function $r$. 
To ensure this, 
one should study the behavior of $\mu_t$ around the boundary $x= 1$
for very small $\gamma$.
Yet, 
when the process $X$ stays close to $1$,
there is a non-trivial competition 
between the amplification 
through the Feynman-Kac penalization 
and the fixation rate at $1$, 
This analysis is beyond the reach of this work.
Even if our conjecture were false, 
the survival of the QSD would mainly rely on a vicinity of $1$,
since~:
\begin{prop}
	\label{TS:prop:sI}
	For any $s, \eps>0$ and measurable bounded $r$, 
	for any $t\ge 1$ sufficiently large~:
	$\quad \PR_{1-\eps} (t < \tau_\eps\bv t < \ext) 
	\rightarrow 1
	\, as\,  \gamma \rightarrow 0,$
	where $\tau_\eps := \inf\{t\ge 0 : X_t \le \eps\}$.
\end{prop} 

Note that in the deterministic limit $(\gamma = 0)$,
there is no more extinction but a convergence to $0$ at exponential rate.
Thus, $\PR_{1-\eps} (t < \tau_{0,1}) \rightarrow 1$ 
as $\gamma  \rightarrow 0$
is to be expected.
\\

The previous conjecture would imply the following result :

"Given any bounded function $r$,
there is a critical value $s_\vee$ such that 
for any $s\ge s_\vee$ and $\gamma \in \R_+$,
$\rho_0\wedge \rho_1 < \rho_\alpha$."

Such result would imply that polymorphism cannot subsist 
when the selection at the individual level is too large.

\section{Discussion on the results on $\mu_t$
	with regards to the individual-based model}
\label{TS:sec:Disc}
\setcounter{eq}{0}

The following discussion
aims at answering to the next questions~: 
\begin{itemize}
	%
	\item Do the asymptotic results of convergence to QSDs 
	and their stability properties
	provide an accurate idea of the basins of attraction ?
	
	\item Is the convergence at exponential rate obtained in the previous convergences
	representative of the observed dynamics ?
	
	\item Are the intra-group fluctuations effectively able 
	to make the inter group selective effects 
	overcome the intra-group selective effects ?
\end{itemize}
In the following discussion, 
we will restrict ourselves to regime of parameters 
for which $(\mu^{m, n}_t)$ may behave similarly as $\mu_t$. 
Referring Theorem \ref{TS:MLim},
we mean that $n$ is thought to be large 
but $m$ to be much larger.
So we look at the regime 
where the inter-group fluctuations 
happen in a much larger time-scale than the intra-group fluctuations. 
Referring to the parameters in the limit $\nu_t$ 
(satisfying the martingale condition \Req{TS:LWF}),
it corresponds to $\gamma_G \ll \gamma_I$.

Subsection  \ref{TS:sec:pureSel} is focused on the two first questions
in the case where the scarcity of intra-group fluctuations
might raise some issues.
We deal with the last one in Subsections \ref{TS:InterG},
supported by  intriguing numerical results
presented in Subsection \ref{TS:Num}.
These conclusions are confronted to more classical results
of weak selection
in  Subsection  \ref{TS:weak_sel}.
Given that these conclusions are quite different,
we elaborate  in  Subsection  \ref{TS:group_sel}
on the specificity of this model
as an instance of the contested notion of "group selection".

An evolutionary perspective is also considered
in Subsection \ref{TS:evol_context},
in a situation where mutations are rare.
We thus imagine a sequence of periods 
in which a single type dominate (the resident)
with brutal transitions
when this type gets replaced
by a mutant subpopulation.
Depending on the mutation effect, 
the pair resident/mutant type could be seen either as C / D or as the opposite.
Of course, such a replacement
is generally rare among the numerous mutations
that are continuously generated,
but we look on a long time-scale where numerous of these events
have happened.
We then try to answer the following question :
Which types of mutations effectively drive the evolutive dynamics of the system ?

Our conjecture is based on the expected stability properties
of the process with a finite yet large population size
around the stable quasi-stationary distributions of the limiting process.
This is compared to the weak-selection regime that is easier to study, 
yet yielding very different conclusions.

Based on these expected stability properties,
we also conjecture 
in Subsection \ref{TS:sec:poly}
that in a regime where a polymorphic QSD
is stable for the dynamics of $\mu A_t$,
it is also a very stable attractor 
for the dynamics in finite yet large population size.
Finally, we  conclude in Subsection \ref{TS:sec:concl}
by the main conclusions of Section \ref{TS:sec:Disc}.

\subsection{Close to the purely selective case.$\quad$}
\label{TS:sec:pureSel}
We first begin the comparison 
by connecting to the results of \cite{LM15}.
We thus focus on the vicinity of their limits, 
namely when $\gamma$ is quite small as compared to $r$ and $s$.
The first prediction that we can get from \cite{LM15}
is that $\mu_t$ goes to $\delta_0$
provided that $\mu_0([1-\eps_\gamma, 1]) = 0$
for a small value of $\eps_\gamma$ 
that we could let tend to $0$ as $\gamma \rightarrow 0$.

\subsubsection{First conjecture and limitations.$\quad$}
Looking at Proposition \ref{TS:prop:sI}
and our conjecture that
$\lim_{\gamma \rightarrow 0} \rho_\alpha(\gamma) 
= \infty$ also holds for any $s>0$ and any bounded function $r$,
we would be in the case $\rho_1< \rho_0< \rho_\alpha$.
From our following Propositions \ref{TS:p:r1r0ra} and \ref{TS:p2:r1r0ra},
one would a priori predict an eventual convergence to $\delta_1$,
with the vast majority of the other groups 
being fixed as pure D groups.
We can however expect 
that the proportion of groups that get fixed
as pure $C$ groups might be very tiny for times
that are not so large.
This would possibly be too tiny 
for the approximation of $(\mu^{m, n}_t)$ by $(\mu_t)$
to be really valid in this vicinity of $1$
for a large range of values of $n$ and $m$. 
Notably, if $m$ is not so large, 
we can not observe in $(\mu^{m, n}_t)$ 
a proportion less than or close to $1/m$ 
since $1/m$ is the contribution of every single group
and the risk of extinction is very large for small sub-populations of groups.

\subsubsection{Simulations with various initial conditions.$\quad$}
Some illustrations obtained by simulations
for this case
are  given in Figure \ref{TS:FD26}-\ref{TS:QSD20}.
\begin{figure}[t]
	\vspace{0.2cm} \hrule
	\vspace{0.2cm}
	\textsl{Close to the purely selective case with $\rho_1 < \rho_0 < \ra$ :}
	\vspace{0.2cm}
	\hrule
	\vspace{0.2cm}
	\begin{minipage}[c]{.48\linewidth}
		\begin{center}
			\includegraphics[width=55mm, height = 27mm]
			{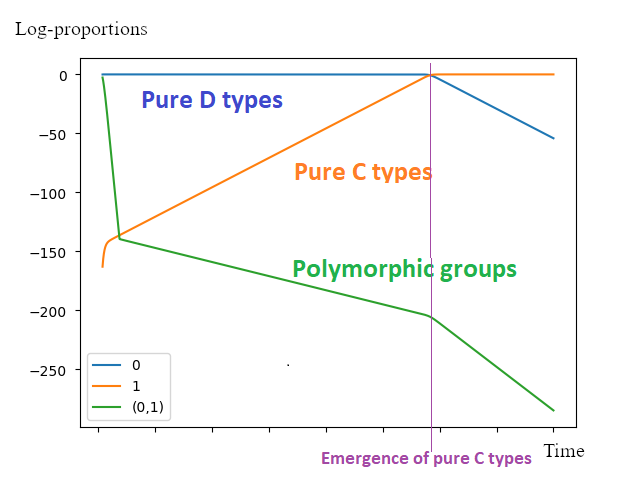}
			\caption{Dynamics of the log-proportions}
			\label{TS:FD26}
		\end{center}
		With the notations of Proposition \ref{TS:DecExt}, 
		%
		the curves presented resp. in blue, orange and green 
		correspond resp. to $\{\log(x^0_t)\}_{t\ge 0}$,
		$\{\log(x^1_t)\}_{t\ge 0}$ and $\{\log(x^\xi_t)\}_{t\ge 0}$.
	\end{minipage}
	\hfill\vrule\hfill
	\begin{minipage}[c]{.48\linewidth}
		\begin{center}
			\includegraphics[width=55mm, height = 27mm]
			{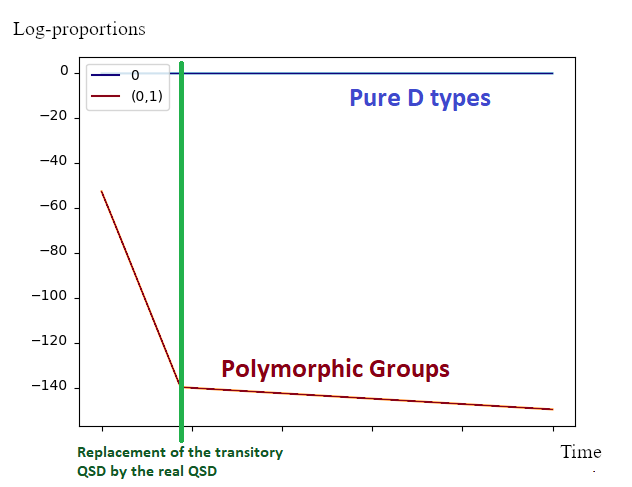}
			\caption{Dynamics of the log-proportions restricted to $[0, 1)$}
			\label{TS:FDL26}
		\end{center}

		The curves presented resp. in blue and red
		correspond resp. to $\{\log(x^0_t)-\log(x^0_t + x^\xi_t)\}_{t\ge 0}$
		and  $\{\log(x^\xi_t)-\log(x^0_t + x^\xi_t)\}_{t\ge 0}$.
	\end{minipage}
	
	\vspace{0.2cm}
	\hrule
	\vspace{0.2cm}

	\begin{minipage}[c]{.46\linewidth}
		\begin{center}
			\includegraphics[width=55mm, height = 27mm]
			{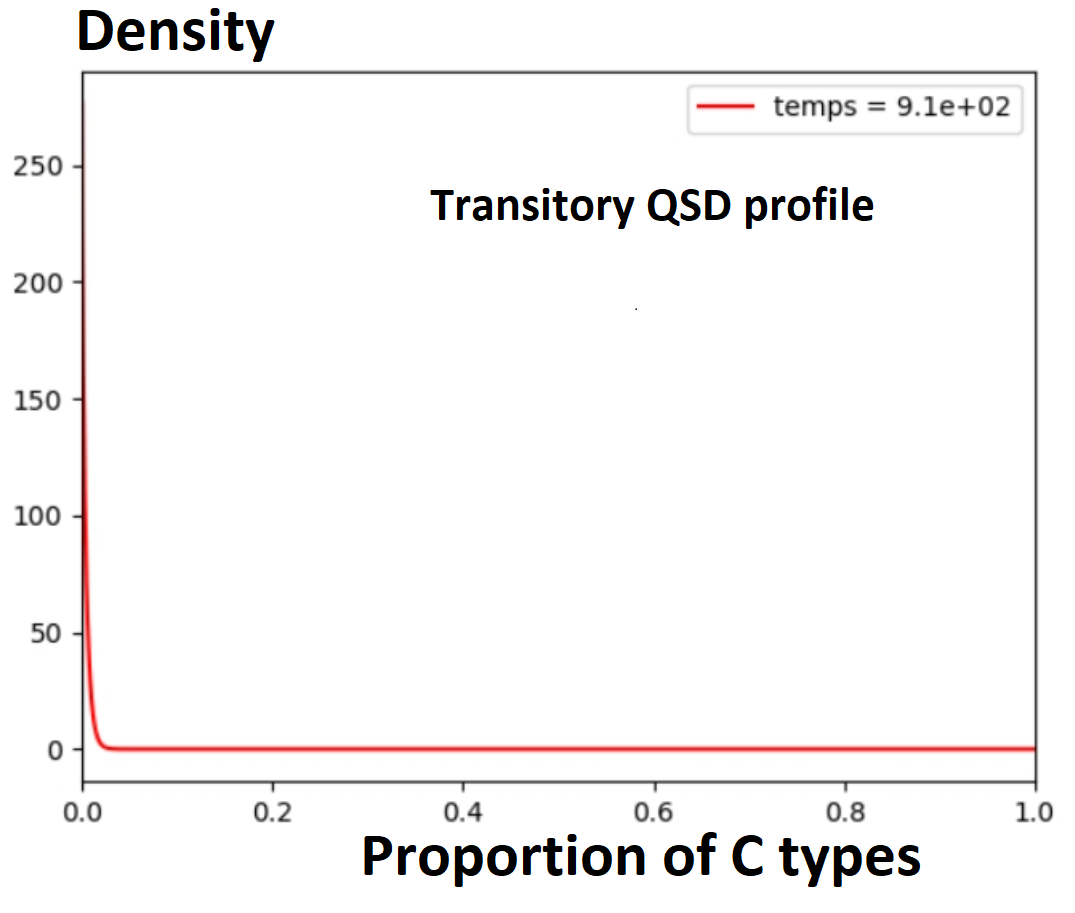}
			\caption{"Transitory QSD" of polymorphic groups}
			\label{TS:tQSD26}
		\end{center}
	\end{minipage}
	\hfill\vrule\hfill
	\begin{minipage}[c]{.46\linewidth}
		\begin{center}
			\includegraphics[width=55mm, height = 27mm]
			{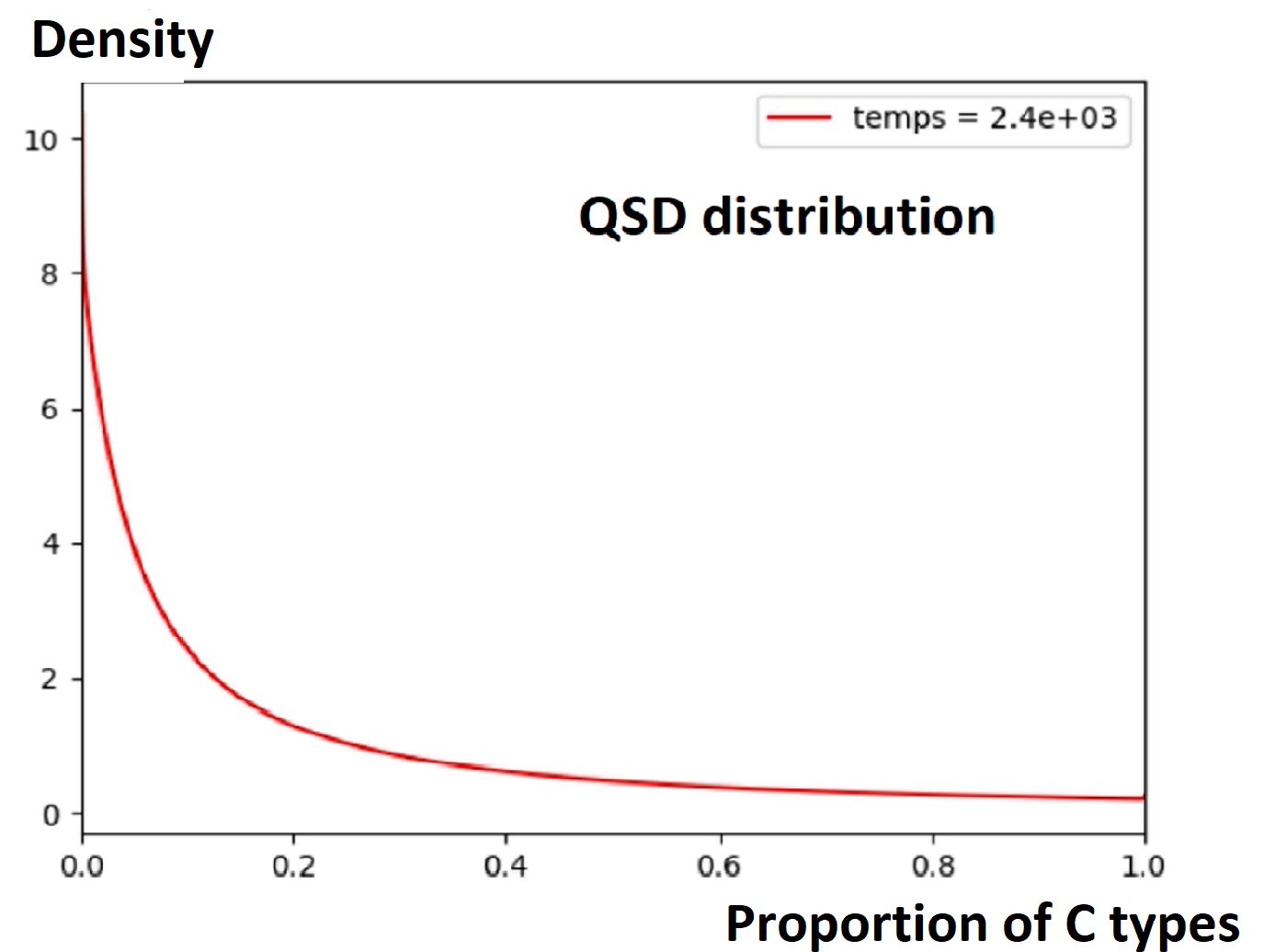}
			\caption{"Final QSD" of polymorphic groups}
			\label{TS:QSD26}
		\end{center}
	\end{minipage}
	\vspace{0.2cm}
	\hrule
	\vspace{0.2cm}
	These results of simulations were obtained with the following parameter values :
	
	$\sqrt{2\gamma} = 0,005$, $s = 0,1$ and $r(x) = 0,005 \ltm x$.
	
	\vspace{0.2cm}\hrule
	\vspace{0.2cm}
\end{figure}
\begin{figure}[t]
	\vspace{0.2cm} \hrule
	\vspace{0.2cm}
	
	\textsl{Close to the purely selective case with $\rho_1 < \ra< \rho_0$  :}
	\vspace{0.2cm}
	\hrule
	\vspace{0.2cm}
	\begin{minipage}[c]{.46\linewidth}
		\begin{center}
			\includegraphics[width=55mm, height = 27mm]
			{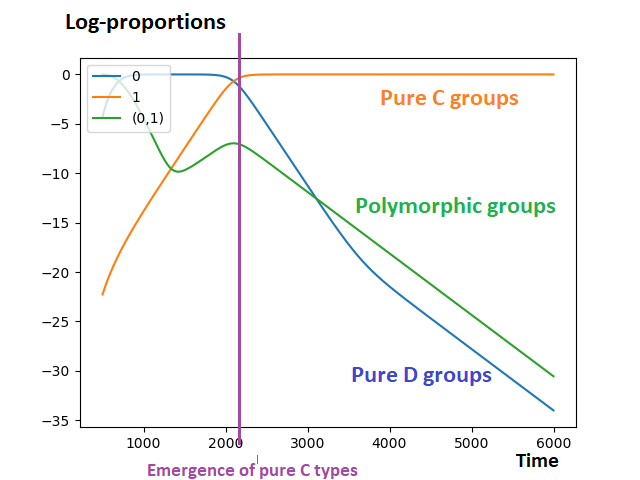}
			\caption{Dynamics of log-proportions}
			\label{TS:FDL20}
		\end{center}	
		The curves presented resp. in blue, orange and green 
		correspond resp. to $\{\log(x^0_t)\}_{t\ge 0}$,
		$\{\log(x^1_t)\}_{t\ge 0}$ and $\{\log(x^\xi_t)\}_{t\ge 0}$.
	\end{minipage}
	\hfill\vrule\hfill
	\begin{minipage}[c]{.46\linewidth}
		\begin{center}
			\includegraphics[width=55mm, height = 27mm]
			{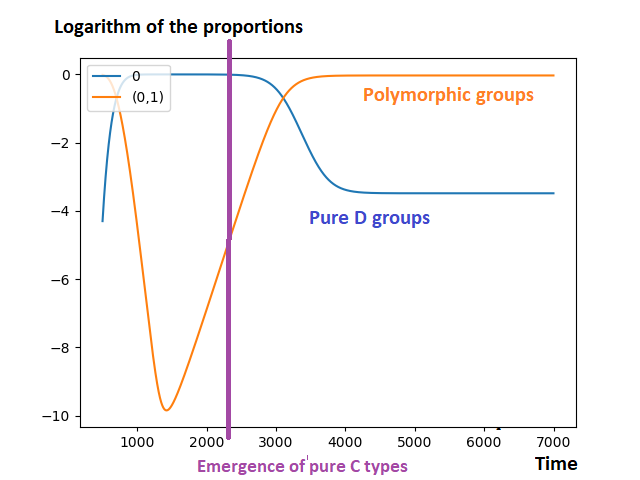}
			\caption{Dynamics of the log-proportions restricted to $[0, 1)$}
			\label{TS:FDLC20}
		\end{center}
		The curves presented resp. in blue and red
		correspond resp. to $\{\log(x^0_t)-\log(x^0_t + x^\xi_t)\}_{t\ge 0}$
		and  $\{\log(x^\xi_t)-\log(x^0_t + x^\xi_t)\}_{t\ge 0}$.	
	\end{minipage}
	
	\vspace{0.2cm}\hrule
	\vspace{0.2cm}
	
	\begin{minipage}[c]{.46\linewidth}
		
		\begin{center}
			\includegraphics[width=55mm, height = 27mm]
			{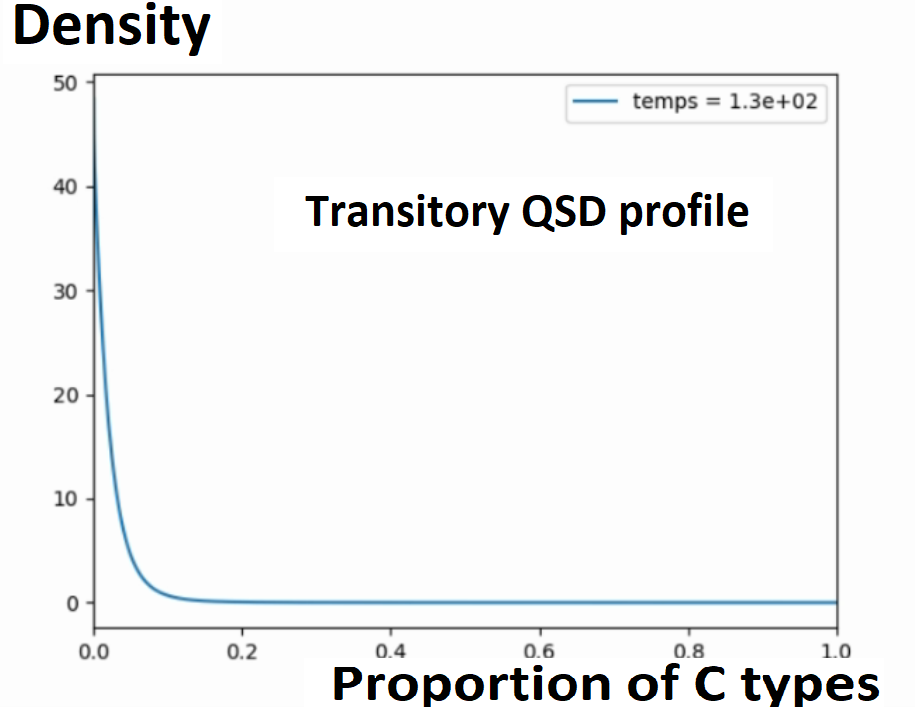}
			\caption{"Transitory QSD" of polymorphic groups}
			\label{TS:tQSD20}
		\end{center}
	\end{minipage}
	\hfill\vrule\hfill
	\begin{minipage}[c]{.46\linewidth}
		\begin{center}
			\includegraphics[width=55mm, height = 27mm]
			{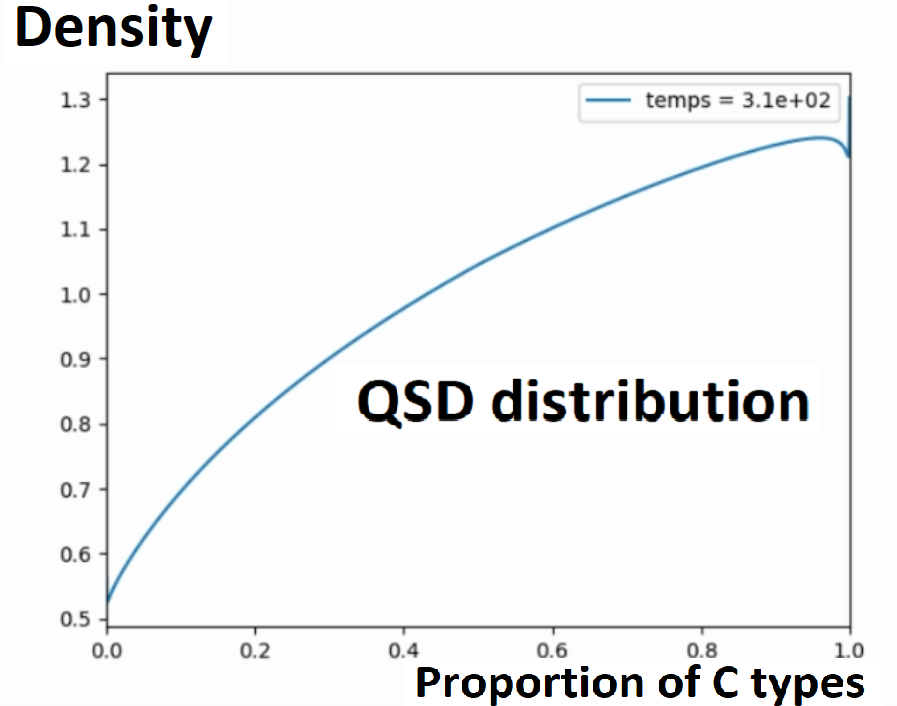}
			\caption{"Final QSD" of polymorphic groups}
			\label{TS:QSD20}
		\end{center}
	\end{minipage}
	\vspace{0.2cm}
	\hrule
	\vspace{0.2cm}
	These results of simulations were obtained with the following parameter values :
	
	$\sqrt{2\gamma} = 0,005$, $s = 0,03$ and $r(x) = 0,1 \ltm x$.
	\vspace{0.2cm} \hrule
	\vspace{0.2cm}
\end{figure}

For this purpose, we tried to derive the dynamics
starting from Dirac initial conditions
with various positions.
Our simulations seem indeed to generally indicate
that for small values of $\gamma$,
and as $t$ goes on,
the distribution $\mu_t$ seems first 
to be attracted by the vicinity of $0$,
where type D individuals prevail.
$\mu_t(0)$ is then very close to 1 
after a time that depends on the initial condition,
see Figure \ref{TS:FD26}.
While looking at the delay in this fixation time 
(say at a proportion $1-\eps$, with a tiny $\eps$)
between different Dirac initial conditions,
we see that it is close to the time needed for the deterministic flow 
to bring the condition furthest to 0 to the closest one.
We also considered the law of  $\mu_0 A^{0,1}_t$,
i.e. the one of $\mu_t$ 
conditioned on polymorphic groups.
It usually stabilizes
for long in some kind of attractive state,
quite concentrated around 0,
cf figure \ref{TS:tQSD26} and \ref{TS:tQSD20},
before the actual QSD emerges.
The mass of this final QSD 
is usually rather around $1$ for not so large $\gamma$
as in \ref{TS:QSD20},
but could present more density around $0$ as in Figure \ref{TS:QSD26}.
However, for small values of $\gamma$,
this profile  usually occur after the emergence and domination
of pure C groups,
and possibly long after !

\subsubsection{The convergence result of $\mu A^{0,1}_t$ to its QSD
	might be of little significance.$\quad$}
As stated in Proposition \ref{TS:p:QSD01},
there is strictly speaking only one QSD for the extinction time $\extOU$.
In the context of the previous paragraph,
this exact QSD seems not to play a significant role 
unless there is initially a non-negligible proportion 
of polymorphic groups with a vast majority 
of C types.
Asymptotically in the diffusion model,
the domination by pure C groups 
seems not to be following the convergence to the QSD
but rather to happen concurrently.
With the notations of Proposition \ref{TS:DecExt},
we mean that
at any time $t$ where the QSD is representative of $\mu A^{1}_t$,
$x^1_t$ is already close to one.
The rate of convergence to $1$ of $x^1_t$
might thus approach $\rho_\alpha - \rho_1$
only at a very late stage.
This stage is even possibly meaningless for finite population sizes
since at some point, there would simply no more polymorphic group.

\subsubsection{An often observed alternative QSD.$\quad$}
We investigate in this paragraph 
the above mentioned concentration effect of $\mu_0 A_t^{0,1}$
to an initial profile close to 0, yet different from the final QSD.
By this, we mean that the profiles presented in Figures \ref{TS:tQSD26} and \ref{TS:tQSD20}
can be observed with little variations
for various initial conditions,
provided the latter has very little mass around 1.
As we can see by comparing those profiles to the final QSD profiles 
in Figures resp. \ref{TS:QSD26} and \ref{TS:QSD20},
they are very different.
The transition from the first ones to the latter attractors
happen quite abruptly after an apparent convergence.
In Figures \ref{TS:FD26}, \ref{TS:FDL26}, \ref{TS:FDL20}, \ref{TS:FDLC20},
we observe a distinct slope for the two regimes of attractions
in the log-proportions resp. $\{\log(x^\xi_t)\}_{t\ge 0}$ 
and $\{\log(x^\xi_t)-\log(x^0_t + x^\xi_t)\}_{t\ge 0}$.
The transitions between the two slopes appear quite brutal
as is the transition between the two attractors.

These elements make us interpret the profiles 
given in Figures \ref{TS:tQSD26} and \ref{TS:tQSD20}
as some sort of alternative QSD,
displaying a transitory phase of convergence at exponential rate 
and characterized by a specific rate of fixation towards $0$.
We presume that this distribution is close to the QSD 
for the extinction event given by the first exit of some interval 
$(0, 1-\eps)$, where the smaller is $\gamma$
and the smaller $\eps$ could be chosen.
The relevance of this interpretation 
relies on the scarcity of transitions towards more cooperative groups
that we examine in the next paragraph.

Note that 
is associated to a large towards pure D groups.
The fixation rate of this "observed QSD" 
provides an estimation of the rate of decay of polymorphic groups
in $\mu_t^{n, m}$ given that it gets absorbed as a population of pure D groups.
Since this fixation rate is a priori large,
the contribution of the deterministic flow towards $0$
is however generally more informative.
This effect of convergence 
should not be of crucial relevance in practice
unless one is interested in the asymptotic rate of decay
when $C$ type individuals are already negligible in their groups.
The fact that this profile emerges
has an analogous implication as the convergence result to $\delta_0$
in Theorem 3 of \cite{LM15} :
both indicate that the fixed state 
with only pure D groups 
is essentially very stable against perturbation.
Notably, introducing groups with a few C type individuals 
would generate hardly anything with high probability.
This will be investigated in Subsection \ref{TS:rob_nsep}
in an evolutionary context.

%


\subsubsection{Scarcity of the transitions towards more cooperation, including towards the real QSD. $\quad$}
The real QSD appears in fact 
at very low density at the beginning.
Since it survives at a better rate 
than the observed distribution,
it simply ends up dominating the whole.
The dynamics is presumably similar 
to the observed one leading to the emergence of pure $C$ groups
as $\log(x^1_t)$ tends almost linearly to 0
in  Figures  \ref{TS:FD26} and \ref{TS:FDL20}.
Transitions from the first effective attractor 
to the real QSD 
can however be shown to be very rare.
This explains why the attraction by the "transitory QSD" 
is so visible early on.

To see how rare these transitions are, 
we use some small value(s) $\eta$
that will serve as a threshold of exceptionality.
For such $\eta$,
we define an analogous version $\bar{\mu}^{\eta}_t$
of the discretized solution of $\mu_t$,
with the main difference 
that we set iteratively to 0
the local elements of mass below $\eta$.
The transition may then be delayed or even unobserved, 
possibly even at very low values of truncation ($10^{-20}$ for instance).
The previously mentioned observed QSD then arises, 
concentrated near the pure D group type.

This method is questionable, 
notably because it depends on the time- and space-discretization.
Yet, with a space-grid of around 100 intervals
and stabilization in less than 10,000 steps,
these exceptional transitions
are exactly what we could alter, 
by truncating parts of them,
so as to exhibit their crucial role.
We also expect that the delay between the untruncated dynamics
and the truncated ones
is due to the fact that the front towards groups enriched in C types
is pulled by these exceptional transitions.
The bulk of more polymorphic groups
is not so much involved
in pushing the proportion of types towards more cooperation.
A more advanced numerical scheme 
would certainly be helpful to quantify the exceptionality 
of the trajectories leading the front.
In this view, this distinction between pulled versus pushed wave
could be observed by introducing neutral markers 
whose density would be followed
as in \cite{RGHK12}.

\subsubsection{As a conclusion, a non-forgettable dependency on the initial conditions.$\quad$}
When $\gamma$ is small 
and there is no group for which $C$ types constitute a large majority,
then, from an ecological point of view,
$D$ individuals have fixed in almost every group.
And even conditionally upon the fact 
that this fixation has not occurred
(or with mutations generating 
new type $C$ individuals),
there is a very stable equilibrium 
with groups dominated by $D$ individuals.

\subsection{The contribution of the intra-group fluctuations for intermediate~$\gamma$.}
\label{TS:InterG}

\subsubsection{A surprising numerical observation}
\label{TS:Num}
\begin{figure}[t]
	\vspace{0.2cm} \hrule
	\vspace{0.2cm}
	\textsl{Illustration of a u-turn :}  
	In these figures, the curves represent the distribution $\mu A^{0,1}_t$ 
	on $(0,1)$
	for different values of $t$.
	By convention, the most recent curve is more strongly apparent,
	so as to create the impression of movement.
	\vspace{0.2cm}
	\hrule
	\vspace{0.2cm}
	\begin{minipage}[c]{.32\linewidth}
		\includegraphics[width=40mm, height = 25mm]
		{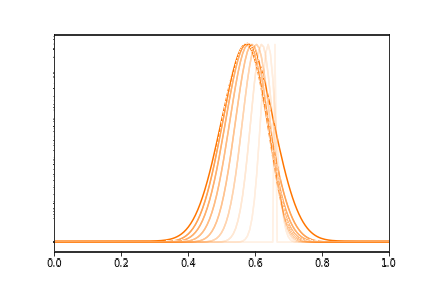}
		\caption{First part of the u-turn with dilatation}
		\label{TS:BuT}
	\end{minipage}
	\hfill\vrule\hfill
	\begin{minipage}[c]{.32\linewidth}
		\includegraphics[width=40mm, height = 25mm]
		{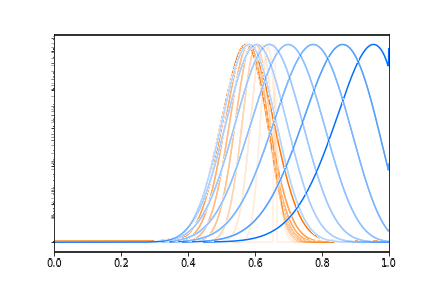}
		\caption{Complete u-turn}
		We used a different color for each of the two parts.
		\label{TS:CuT}
	\end{minipage}
	\hfill\vrule\hfill
	\begin{minipage}[c]{.32\linewidth}
		\includegraphics[width=40mm, height = 25mm]
		{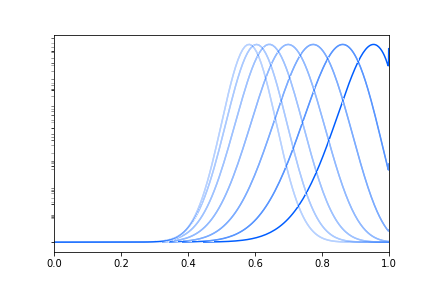}
		\caption{Second part of the u-turn}
		\label{TS:AuT}
	\end{minipage}
	\vspace{0.2cm}
	\hrule
	\vspace{0.2cm}
	These results of simulations were obtained with the following parameter values :
	
	$\sqrt{2\gamma} = 0,02$, $s = 0,1$ and $r(x) = 0,5 \ltm x$.
	\vspace{0.2cm} \hrule
	\vspace{0.2cm}
\end{figure}

There is another interesting behavior 
for not so small values of $\gamma$,
a linear growth rate $r$ (increasing)
and a Dirac mass in the middle of the interval
as initial condition.
At the beginning of its dynamics,
$\mu_t$ is close to a Gaussian distribution 
with an expanding variance and a drift.
Except that the variance is larger,
it seems first to behave as in the case of small $\gamma$
and for initial conditions rather close to $0$, 
not much difference can be observed.
Yet, although the drift is always first directed towards $0$,
we may see a u-turn after a while,
with a drift now seemingly directed towards $1$.
It creates the impression
that the drift at the individual level is changing,
while it is in fact the selective effects at group level 
that starts to play a significant role.
The more diverse the distribution is,
the more these effects differentiate between these different realizations
and the larger is this additional drift. 

\begin{figure}[t]
	\vspace{0.2cm} \hrule
	\vspace{0.2cm}
	\textsl{Significant intra-group fluctuations with $\rho_1 < \ra < \rho_0$ :}         	
	\vspace{0.2cm}
	\hrule
	\vspace{0.2cm}
	
	\begin{minipage}[c]{.46\linewidth}
		\begin{center}
			\includegraphics[width=55mm, height = 29mm]
			{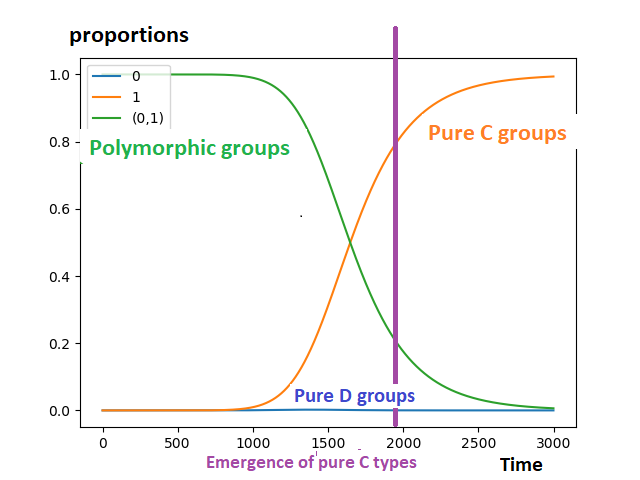}
			\caption{Dynamics of the proportions}
			\label{TS:FD21}
		\end{center}
	\end{minipage}
	\hfill\vrule\hfill
	\begin{minipage}[c]{.46\linewidth}
		\begin{center}
			\includegraphics[width=55mm, height = 29mm]
			{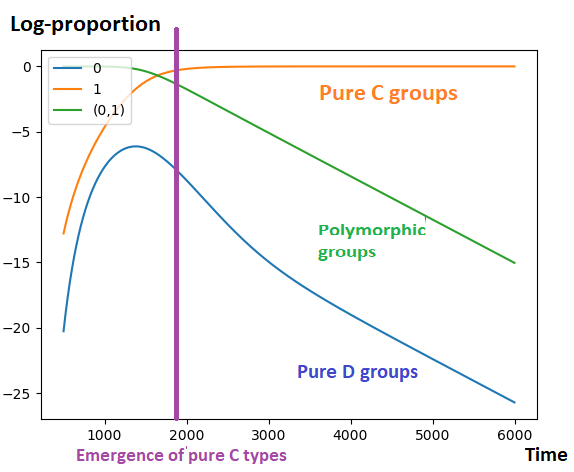}
			\caption{Dynamics of the log-proportions}
			\label{TS:FDL21}
		\end{center}
	\end{minipage}
	
	\vspace{0.2cm} \hrule
	\vspace{0.2cm}
	
	\begin{minipage}[c]{.48\linewidth}
		\begin{center}
			\includegraphics[width=55mm, height = 29mm]
			{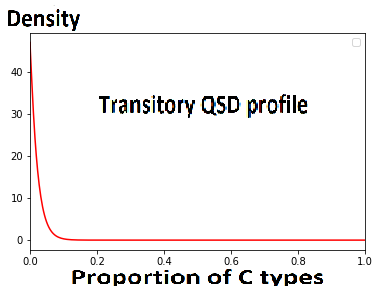}
			\caption{"Transitory QSD"}
			\label{TS:tQSD21}
		\end{center}
	\end{minipage}
	\hfill\vrule\hfill
	\begin{minipage}[c]{.48\linewidth}
		\begin{center}
			\includegraphics[width=55mm, height = 29mm]
			{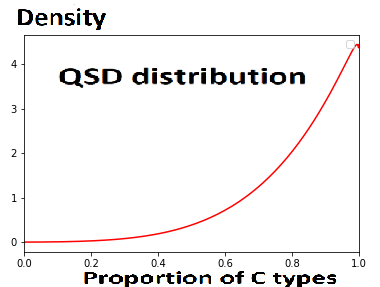}
			\caption{Final QSD}
			\label{TS:QSD21}
		\end{center}
	\end{minipage}
	\vspace{0.2cm}
	\hrule
	\vspace{0.2cm}
	These results of simulations were obtained with the following parameter values :
	
	$\sqrt{2\gamma} = 0,02$, $s = 0,1$ and $r(x) = 0,1 \ltm x$.
	\vspace{0.2cm} \hrule
	
\end{figure}

In such a case, 
the role of the real QSD can be much more significant
as can be seen in Figures \ref{TS:FD21}-\ref{TS:FDL21}.
The convergence to $1$ 
is much more robust against the truncations 
(that we implement as mentioned in the previous paragraph).
We thus expect the estimation of the typical dynamics of $\mu^{m, n}$
by the one of $\mu_t$ to be relatively accurate in this context.
In a future work, 
we plan to validate  
by simulations
notably this assertion
and to be more specific regarding the associated range of parameters
(including the initial condition).
This would demonstrate the crucial importance
of having sufficiently large intra-group fluctuations
for selective effects at group level to be significant.
%
%
\subsubsection{Weak selective effects : $r$ and $s$ small as compared to $\gamma$.$\quad$}
\label{TS:weak_sel}

For simplicity, 
consider the case where $r$ is linear :
$r(x) := r_1\, x,$ for $r_1>0$.
We then may think of the selective effect
as type $C$ individuals distributing 
a reproduction benefit of $r_1 / n$
to all the $n$ individuals in their group.
The specificity is however that 
the whole group is to be duplicated at once
at the reproduction event at group level.
Since we forbid any transmission between groups,
one could expect a relatedness of $1$ 
between any two individuals of the same group,
and $0$ between individuals of different groups.
Under weak selection, 
this may lead to the prediction 
that the mutation is positively selected
provided that $r_1 > s$
(in agreement with eq. 1 in \cite{TN06}, 
with a much larger number of groups than the size of each group,
since only inter-group fluctuations are kept).
Since in our model, 
we assumed that the random fluctuations in group reproduction
are negligible,
we cannot really approach the case of weak selection. 

Such weak mutation assumption
corresponds nonetheless also to $\gamma \rightarrow \infty$.
So what happens in practice in our model in this limit ?
Note first that this implies a separation 
in the time scales of fixation 
inside one group and among groups.
Such time-scale separation is classically assumed 
in the context of  weak selection, recall e.g. \cite{O08}.
In the time scale of interactions between groups,
we can assume with almost no restriction that all the groups have fixed,
so that only remains the competition between pure groups.
Indeed, 
all the individual in the group 
are then strongly related
(with a very close common ancestor).
But the crucial parameter of interest for predicting the outcome
is the initial proportion of pure C groups 
that can be established early on. 
If this is non-negligible, 
pure C groups will very probably prevail 
whatever the respective values of $r$ and $s$ are.

\subsubsection{The contested notion of "group selection".$\quad$}
\label{TS:group_sel}
For such a model 
where the selective effects 
are clearly associated with two hierarchical levels
without transmission between groups,
the notion of selective effects at group level
does not appear so ambiguous.
It is closely related to the more common expression of "group selection",
except that the interpretation of the latter has been quite diverse
depending on the authors.
Notably, we have to mention 
the recurrent discussions on the confusions 
brought about by this notion 
(see \cite{W+07}, \cite{W+08}, \cite{W+11}, \cite{L+07}).

Notably, it has been  argued 
that one shall rather 
relate to the kin selection formalism,
which means that one shall weight 
the selective advantage 
that Cooperators give to other individuals
by the relatedness they share with the recipient.
On the other hand,
our formalism focus on the group level
by specifying the dynamics of a randomly chosen group.
We are convinced that a formalism in terms of relatedness
could produce an equivalent description
but we have not found it particularly helpful.
We would happily welcome any suggestion 
in this regard,
as we find it very intricate in our model.
Notably,
the notion of relatedness 
is here much more difficult to describe,
even in the simplification of weak selection.
This is strongly linked to the fact that individuals
cannot transfer from one group to another
while the internal selective advantage
is relative to the mean selective advantage 
inside the group.
In addition, 
fixation events inside groups 
are a specificity of this model.
The difficulty is also to be expected given that
we do not assume,
contrary to the usual considerations 
on the direction of selection,
that the selective effects 
are on different time-scales
(see e.g.  \cite{O08})
nor weak
(see e.g.  \cite{ML14}, \cite{TN06}).
Already without such population structure, 
different selection intensities may lead 
to qualitatively different relations of dominance, 
see for instance \cite{ML14}, \cite{WGHT13}, \cite{WD19}.

\subsection{Evolutionary context.$\quad$}
\label{TS:evol_context}

\subsubsection{Context and motivations}
We are also interested in looking at an evolutionary scenario
where mutations that are well-separated in time 
have each some probability to invade a population with a given type.
The invasion then corresponds 
to a replacement of the resident type
as the successful mutant type.
In this context,
one has first to choose if $C$ or $D$ plays the role of the resident type.
Let us first say it is $D$ type.
We thus consider the case where only one individual is of type $C$
and compare the probability that an homogeneous population of type $C$ emerges
to the neutral case 
where that probability would be $1/(m\ltm n)$.
We also consider the reverse case 
where a single mutant of type $D$ invade a resident population of type $C$.
Then, we can compare both to infer both the direction favored by natural selection,
the time-scale of this evolution and whether the main contributions shall come 
from mutations of small or large effects.

In the following,
we need to keep a quantification of the fluctuations at group level,
with the parameter $\gamma_G>0$ (cf the model given by equation \Req{TS:piT}).
Recall that $\gamma$ is to some extend its equivalent at individual level.
Since our limiting model correspond to vanishing $\gamma_G$, 
we shall restrict our analysis to the cases where $\gamma_G << \gamma$
(corresponding to much larger number of groups than of individuals inside groups, 
i.e. $m\gg n$).
Note that these values are expected to stay fixed for a given population,
even while undergoing several events of fixation.
We shall first consider the case were we can justify a separation of time-scales,
that makes the study much easier.
We then investigate in Subsection \ref{TS:rob_nsep}
to which extend these results generalize 
when selective effects and intra-group fluctuations 
occur concurrently.

\subsubsection{Separated time-scales}

Note that in this neutral case, 
the invasion probability is respectively 
of order $1/n$ and $1/m$ for the fixation within the group
and of one group among the population of groups.
The invasion probability of one individual into the whole population
is merely the product of these two events. 
It is generally due to the fact that at some large time, 
every individual in the population
will be issued from a unique common ancestor at generation $0$,
which under neutrality is clearly taken uniformly at random.
Such decomposition by product is simpler to analyze and well justified 
in a context where both events happen on separate time-scales.
This is why we first focus on this case. 
We shall then pursue the analysis for more intricate cases 
where the selective effects at group level act on the same time-scale 
as the internal dynamics of a group.

\paragraph{Estimation of the probability of invasions at individual level.$\quad$}
Because there is this separation of time-scales,
one can exploit classical results 
for estimating the probability of invasion
at both individual and group levels.
At least for the selection within group,
we may exploit the explicit formula for probability of invasion
starting from a proportion $x\in (0,1)$ 
in the solution $(X_t)$ of the limiting equation \Req{TS:Xdef}.
The formula, first obtained by Malécot 
and often presented in the context of Kimura's diffusive approximation, 
takes the form :
\begin{equation}\EQn{TS:KimFix}{}
	\PR_x(\tau_1  < \tau_0) = \dfrac{\exp[sx/\gamma]-1}{\exp[s/\gamma]-1}.
\end{equation}
This can be obtained by identifying it as the only solution 
of $\cL_{WF} u = 0$ 
with boundary conditions $u(x) = 0$ and $u(1) = 1$.
This expression is then equivalent as $x$ tends to $0$ to 
$(s/\gamma)\ltm (e^{s/\gamma} -1)^{-1}\ltm x$.
Although a precise justification would require 
a careful analysis of the process when the $C$ sub-population is still negligible,
we can expect that selective effects do not play a consequent role 
in this first step.
We may thus expect a fixation probability 
starting from only one type $C$ individuals to be well-approximated by
\begin{equation}\EQn{TS:pIdc}{}
	\pi^I_{D\mapsto C} = \frac{s/\gamma}{n\, (e^{s/\gamma} -1)}
	\approx (1/n)\ltm (1 - \frac{s}{2 \gamma}) 
	\with s/\gamma \ll 1.
\end{equation}
Note that except in this last step of approximation, 
we only exploited the approximation by $X$ for $n$ large
and the separation of time-scale without restriction on $s$ or $\gamma$.
By changing $s$ into $-s$, we would obtain the probability of invasion
of $C$ type residents by $D$ type mutants.

\paragraph{Estimations of the probability of invasions under a weak selection assumption.$\quad$}

The same reasoning can then be applied for the subsequent fixation 
of the pure $C$ group in a resident population of pure $D$ groups,
leading to an overall invasion probability well-approximated by :
\begin{align*}\EQn{TS:piDC}{\pi_{D\mapsto C}}
		\pi_{D\mapsto C} 
		&= \pi^I_{D\mapsto C} \ltm  \pi^G_{D\mapsto C}
		= \frac{s/\gamma}{n\, (e^{s/\gamma} -1)}
		\ltm \frac{(\rho_0 -\rho_1)/\gamma_G}{m\, (e^{(\rho_0 -\rho_1)/\gamma_G} -1)}
		\\&
		= \frac{1}{n\ltm m} 
		\times \frac{s/\gamma}{e^{s/\gamma} -1}
		\ltm \frac{r_1/\gamma_G}{1-e^{-r_1/\gamma_G}}
		\approx \frac{1}{n\ltm m}  \ltm (1 + (1/2)\ltm (r_1/\gamma_G - s/\gamma)),
		\\
		\pi_{C\mapsto D} 
		&= \frac{1}{n\ltm m} 
		\times \frac{s/\gamma}{1-e^{-s/\gamma}}
		\ltm \frac{r_1/\gamma_G}{e^{r_1/\gamma_G}-1}
		\approx \frac{1}{n\ltm m}  \ltm (1 - (1/2)\ltm (r_1/\gamma_G - s/\gamma)),
\end{align*}
where the last approximations assume $(r_1/\gamma_G)\vee (s/\gamma) \ll 1$.
Assuming that mutations in both directions happen with the same law and intensity,
well-separated and with small selective effects,
we would then expect an evolutionary drift in the direction of cooperation
provided $r_1/\gamma_G > s/\gamma$.

\paragraph{A specific criterion for the direction of selection.$\quad$}
We see that the above-mentioned prediction 
is not the one we find here, 
with an additional implication 
of the levels of fluctuations $\gamma$ and $\gamma_G$.
Remark that the inverse of these quantities is usually referred to 
as the effective population size.
This can be derived from Theorem \ref{TS:MLim}
by recalling  that $\gamma \approx \bar{\gamma} / n$
with $\bar{\gamma}$ the actual reproduction rate of the individuals.
Assuming $\bar{\gamma}$ of order $1$ 
then implies that $\gamma$ scales as the inverse of the group size.
A similar interpretation scales $\gamma_G$ as the inverse of $m$. 
Besides, 
we recall that we justified this separation of time-scales by assuming $\gamma\gg \gamma_G$.
So  our condition $r_1/\gamma_G > s/\gamma$
is in this context very different from the naive one $r_1>s$,
obtained by comparing the benefits in terms of reproduction rate
(see the above paragraph on weak selective effects).
We see here 
that even small selective effects at group level
can effectively outcompete much larger selective effects at individual level,
provided the number of groups is sufficiently large
as compared to the number of individuals inside groups.

\paragraph{Separation of time-scales with strong selection}

The product expression in \Req{TS:piDC} 
can be exploited as long as the separation of time-scales occur.
Because each level is then treated separately, 
the fate of the system almost exclusively depend on the ratio $r_1/\gamma_G$ and $s/\gamma$,
that we now allow to be non negligible.
\\

Naturally, a type 
has a very small probability of being invaded by a mutant subpopulation
provided the type is clearly advantageous at group level 
($r/\gamma_G$ sufficiently large)
while only weakly deleterious at individual level 
($s/\gamma$ sufficiently small). 
Indeed, we find in such a case :
\begin{equation}\EQn{TS:piGCD}{}
	\pi^G_{C\mapsto D} 
	\approx (1/m) \times (r_1/\gamma_G)\ltm e^{-r_1/\gamma_G} \ll (1/m).
\end{equation}
while $\pi^I_{C\mapsto D}$ is not so far from $1/n$.

Similarly, such an advantage at group level 
is also more likely to fix, although the effect is less strong as $r_1/\gamma_G$
gets large. We find indeed :
\begin{equation}\EQn{TS:piGDC}{}
	\pi^G_{D\mapsto C} 
	\approx (1/m) \times (r_1/\gamma_G),
\end{equation}
while  $\pi^I_{D\mapsto C}$ is again not so far from $1/n$.
\\

Note that in such regime of separation of time-scales, 
the individual and group levels play symmetric roles, 
so that one can easily reverse the argument for strong selective effects
at individual level against weak effects at group level.

By symmetry and assuming that the mutations
is as likely to give an advantage $s$ at individual level 
with a penalty $r$ at group  level
as to give $s$ as a penalty and $r$ as an advantage,
we find again that the direction of selection is given by the higher ratio between $r/\gamma_G$ and $s/\gamma$.
Under the condition of a constant ratio $(r\ltm \gamma / s\ltm \gamma_I)$,
the speed of selection appears to vanish as the selective effects get larger and larger.
Moreover, assuming for instance $r/\gamma_G > s/\gamma$, 
the fixation probability of a mutation deleterious at group level 
gets more and more negligible as compared to a mutation deleterious at individual level :
as for a ratchet, 
selection hardly let any change of its direction happen.
From the previous calculations, 
we indeed see that the probability of seeing no fixation
at some level
gets much smaller
than the probability of fixation increases at the other.
\\

%
%
%
%

\subsubsection{This robustness to invasions shall extend 
	even without the separation of time-scales.
	$\quad$}
\label{TS:rob_nsep}

The situation is much more complicated to analyze in this case 
where $\gamma$ is not so large,
so that both selective effects are competing simultaneously.
Nonetheless, 
in the limiting process we described,
as soon as there is a  non-negligible proportion
of pure  $C$ groups 
in the population,
$D$ individuals 
simply cannot completely replace $C$ type individuals.
The complete invasion by $D$ type individuals 
is impossible,
even in the case where $0< r_1 \ll s$ (nor $r_1\ll (s/\gamma)$).
With a large yet finite number $m$ of groups,
we expect that it would be possible 
to interpret such invasion as a large deviation result
of the process $\nu$ (see section 1).
Referring to classical literature 
on the subject of large deviations, 
notably Section 5 of \cite{DZ98},
the associated probability is thus likely 
to be exponentially small with increasing $m$.
At least, this rate of decay
is what we have obtained in \Req{TS:piGCD}, 
recalling $m = O(1/\gamma_G)$.
So we conjecture that this strong resistance to invasion 
by Defectors is very general as long as $\gamma_G$ is sufficiently small,
that is $m$ sufficiently large.
The selective effects at individual level 
shall certainly make the invasion more frequent 
than it would be without it, 
yet generally not sufficiently to compensate for this exponential reduction.

\paragraph{The invasion by $C$ type mutants is not as strongly selected 
	against when $\gamma$ is not so large.$\quad$}
When $\gamma$ is sufficiently large to keep non-negligible
the probability of invasion of its group 
by a mutant C type individual,
we shall retrieve the ratchet effect :
neglecting the effects of weak mutations,
invasions by Cooperative individuals 
should drive the dynamics of selection
towards more cooperation,
while invasions by Defectors scarcely occur.
We would have the same effect if mutations towards more cooperation
were strong and simply more advantageous in an homogeneous population (without group structure).
Globally, the dynamics is driven also by these weak mutations
and the contribution of both weak and strong selection
a priori depends on the specific situation of study.

\paragraph{Robustness to invasion by $C$-type mutants when $\gamma$ is small.$\quad$}
On the contrary, 
we also have a similar resistance 
when $\gamma$ is sufficiently small,
as noted in our Section \ref{TS:sec:pureSel}.
Even if $r_1 \gg s >0$, 
the invasion of some group
with mainly type $D$ individuals
by some type $C$ individuals
relies on so exceptional events 
that it seems biologically almost impossible.
We expect that it would be possible 
to interpret such invasion as a large deviation result
of the process $X$, cf \Req{TS:Xdef}.
In practice, we thus predict
that the probability of such invasion 
shall decrease exponentially with $\gamma = O(1/n)$.

The case where $r_1$ is of order $s/\gamma$ 
is possibly more intricate and would require further consideration.
It might lead to a specific optimization problem as described 
in \cite{Henry},
where the cost of deviating $X$ has to be balanced 
with the amplification through $r$.
The shape of the function $r$ 
would then play a much  more significant role.

\paragraph{A general robustness to invasion when $\gamma$ is small.$\quad$}
Considering both directions of invasion, 
it seems that strong selective effects 
are strongly selected against
in very large populations,
whatever the level of selection they favor,
as long as they are detrimental 
for one level of selection. 
Indeed, 
both invasion probabilities scale 
as exponentially small in the population sizes 
(of the detrimental level of selection),
which is much more stringent than the order $O(1/(n\times m))$
of nearly neutral mutations.
We expect it to extend with possibly 
more levels of "selection".
Thus,
beside the effect of mutations favoring its carriers at both levels 
(but not necessarily equivalently),
the trade-off between selective effects 
at different levels shall be 
driven mostly
by weak selective effects.

Note that these conclusions seem quite robust to more general forms 
of functions $r$,
provided $\rho_1< \rho_0 < \rho_\alpha$
with a QSD $\alpha$ very concentrated around $1$
(for robustness against invasions by D types)
and provided there is an apparent QSD near 0 
(for robustness against invasions by C types).
Based on our first simulation results,
this effect of concentration seems to be robust
as long as $r$ is increasing with increasing cooperation.
This suggests also that 
it is not so crucial 
that $1$ is an absorbing state.
Thus, 
the above conclusions
shall be maintained even with some  transmission between groups,
as long as their rate is sufficiently small.

\paragraph{A similar robustness for polymorphic QSD.$\quad$}
\label{TS:sec:poly}
Similarly, the fixation of a polymorphic population 
with profile $\alpha$ such that $\rho_\alpha < \rho_0\wedge \rho_1$
is likely to be an exceptional event
as compared to the time-scale at which the transitory profile evolves.
We mean that the profile of $\mu^{m, n}$ shall remain very close to $\alpha$,
with a very quick regulation of the random perturbations
when $m$ is large.
Our confidence originates again from the comparison
with Large Deviations results.
Even the events of fixation, 
that are possibly much less negligible for finite $m$,
are not expected to be significant. 
The pure groups have a lower progeny 
and do not contribute much to the dynamics on the long term.
%
%
By construction, 
this case corresponds to $r$ being a function 
with at least a strict maximum
inside $(0,1)$. 
Selective effects must favor polymorphism directly 
and not only conflict with selective effects at individual level.
Otherwise, $\rho_\alpha < \rho_1$ is excluded by Lemma \ref{TS:ExtRho}.
We could also imagine more general selective effects at individual level 
with a frequency-dependency. 
This would possibly also entail a stable polymorphic QSD.

Note that we also assume here that $\gamma$ is not negligible.
Again, referring to \cite{Henry},
it is possible that 
the description gets much more tricky 
in the case where $\gamma$ is small 
but $s$ also so that the function $r$ scales as $s/\gamma$.

\paragraph{Difficulties in relating to kin selection.$\quad$}
Considering more general $r$, 
notably with a maximal value in the interior $(0,1)$,
the approach of kin selection becomes even less clear.
It seems required to deal with another definition of relatedness,
like the one given in \cite{G06},
with much more complexity.

\subsection{Main conclusion of Subsection \ref{TS:sec:Disc}}
\label{TS:sec:concl}
\IEn
\ENm  In any case where $\delta_1$ is stable, 
for an initial condition
with enough highly cooperative groups
so that the diffusion in them
generates a non-negligible proportion of pure $C$ groups,
these groups do eventually invade the population with a high probability.
The range of parameters that induce this effect 
strongly depends on the initial condition if cooperation is initially well-established.

\ENm If intra-group fluctuations of population size are small,
a very large proportion of the groups become visibly
increasingly dominated by $D$ individuals,
even though the pure $D$ groups are the worst at reproducing.
Introducing in the model rare mutations from individuals $D$ to individuals $C$ 
probably wouldn't make much of a difference.

\ENm The case $\rho_1 < \ra <\rho_0$ characterizes the fact 
that a subpopulation of polymorphic groups 
is able to maintain itself better than pure $D$ groups.
This means that domination by $D$ individuals would still be prevented,
even if we disrupt the model 
by introducing rare events of migration between groups
or mutations from $C$ individuals to $D$ individuals.
This should hold true as long as the polymorphic QSD $\alpha$ is rapidly approached.

\ENm Likewise, in the case of $\ra< \rho_1< \rho_0$,
the polymorphic state is a priori very stable and
$\mu_t$ tends to him 
provided that the proportion of polymorphic groups
was not violently reduced from the very beginning.

\ENm
However, it is unclear 
that the area of attraction of this quasi-stationary distribution
is far beyond the vicinity of the purely cooperative state.
Convergence towards it and its influence has potentially no significant effect
if the proportion of the groups with a majority of $C$ individuals is too small.
A transient attractor is then likely to appear for $\mu_t$
restricted to $(0,1)$, widely supported on a neighborhood of the purely $D$ state.
It can be interpreted in practice as an alternative quasi-stationary distribution,
with an extinction rate generally higher than $\rho_0$.
This happens when the transition
from a group dominated by $D$ individuals to one with 
mostly $C$ type individuals
is too costly as a deviation from the process $X$.
The difference in extinction rates
is then not able to quickly offset
the very small probability of such a transition.
Transitions in the vicinity of the less stable "pseudo"-QSD 
to the "real" QSD (more stable and around~1)
take a non-negligible time
so that the attraction to the pseudo-QSD is clearly visible.

\ENm In an evolutionary perspective, 
for very strong selection effects
as compared to genetic fluctuations,
the status quo situation is more likely to prevail:
the probability of fixing a mutation 
which puts a burden on its holder at some selection level
is only slightly compensated by the advantage 
that this mutation could bring at another level. 

\ENm Still from an evolutionary point of view,
cooperation is favored
by any increase of the genetic fluctuations at intra-group level.
It namely corresponds to a greater kinship between individuals in the same group.
Assume that the genetic fluctuations bewteen groups are kept small,
corresponding to large values of $m$,
while the genetic fluctuations at intra-group level are non-negligible.
Then, mutations that put the group of those who carry them at a disadvantage
have a much lower probability of fixation
as compared to the mutations that put those who carry them at a disadvantage inside their group:
the relation of comparison between the strength of these effects 
in the weak selection context
seems to lose its relevance here
when the effects combine.

\section{Proof of the results of Section \ref{TS:sec:QSD}}
\label{TS:sec:PfQSD}
\setcounter{eq}{0}

\subsection{Proof of Proposition \ref{TS:p:QSD01} :
	characterization of $\alpha$ on $(0,1)$}

We rely on the method used in \cite{AV_QSD} 
and more precisely on the proof of 
the second illustration
presented in Subsection 4.2 of \cite{AV_QSD} 
to ensure that~:
\begin{align*}
		\Exq{\chi>0}
		\frl{n\ge 1, \xi>0} 
		\Exq{C_{n,\xi}>0}
		\frlq{\mu \in \MnxTS^{0,1}}
		\NTV{\mu A^{01}_t
			- \alpha}
		\le C_{n,\xi} \, \exp[-\chi\, t].
		\EQn{TS:CVMnx}{}
\end{align*}
The diffusion is indeed regular on any 
$\cD_n := [1/n, 1-1/n]$ (for $n\ge 3$) 
so that applying the Harnack inequality,
we prove similarly as in \cite{AV_QSD} 
that for any choice of $0< t_M < t_c$, 
there exists $c_M>0$ such that for any $x\in \cD_n$ :
\begin{align*}
		&\PR_x\lp X_{t_{M}} \in dx\pv
		t_{M} < \ext^{n+1} \rp
		\ge c_{M}\, \PR_{1/2}\lp X_{t_c} \in dx\bv
		t_c < \ext^{3} \rp 
		:= c_{M}\,\zeta(dx),
		\EQn{TS:alc}{}
		\\& \with 
		\ext^n := \inf\Lbr t >0 \bv X_t \notin \cD_n\Rbr.
\end{align*}

We refer to the step 4 of the proof given in Sect. 4 of \cite{ChpLyap}
to ensure that for any $n\ge 3$ and $t>0$,
there exists $c_n>0$ such that ~:
\begin{align*}
		\frlq{x, y \in \cD_n}
		\PR_x\lp X_t \in dx\pv
		t < \extOU\rp
		\le c_n\, \PR_y\lp X_t \in dx\pv
		t < \extOU\rp.
		\EQn{TS:AH}{}
\end{align*}

Next, we prove that the process $X$ cannot maintain itself close to the boundary :
\begin{lem}\label{TS:lem:esc}
	For any $\rho>0$, there exists $E = \cD_{n_{E}}$ such that ~:
	\begin{align*}
			\sup_{x\in (0,1)} \E_x \exp[ \rho V_{E}] < \infty
			\where  V_{E} :=  \extOU \wedge \inf\Lbr t>0\pv X_t \in E \Rbr.
			\EQn{TS:etZ}{}
	\end{align*}
\end{lem}

Applying Theorem 2.1 in  \cite{AV_QSD} with \Req{TS:alc}, \Req{TS:AH} and \Req{TS:etZ},
noting also that condition (A0) on $\Lbr\cD_n\Rbr$ is clearly satisfied,
concludes the proof of \Req{TS:CVMnx}.
The results on the capacity of survival comes from Theorem 2.2.
\\

To end the proof of Proposition \ref{TS:p:QSD01}, 
it is sufficient to ensure the following lemma
\begin{lem}
	\label{TS:esc1}
	There exists $n_{B}\ge3$, $\xi_{B}, t_{B}>0$ such that ~:
	\begin{align*}
			\frlq{\mu \in \MoneM} 
			\mu A^{01}_{t_{B}} \in \M_{n_{B}, \xi_{B}}.
	\end{align*}
\end{lem}

This can be done exactly as in step 1, Section 5.1 of \cite{Ch1D},
by handling precisely with 
the vicinities of $0$ and of $1$.


\paragraph{Proof of Lemma \ref{TS:lem:esc}}

The core of the proof is the well-known fact
that for any $t>0$,
\begin{align*}
		\PR_x (t< \tau_{0,1}) \cvz{x} 0\mVg
		\quad \PR_x (t< \tau_{0,1})\rightarrow 0
		\text{ as } x\rightarrow 1
		\EQn{TS:ExtV}{}
\end{align*}
(see notably Theorem 3.4 and 3.7 in \cite{Ch1D}
for a much more precise estimate 
of the extinction on the boundaries).

Let $\rho>0$.
We fix then arbitrarily $t = 1$ 
and deduce from \Req{TS:ExtV}
that for $n_E$ sufficiently large :
\begin{equation*}
	\frlq{x\in (0, 1)}
	\PR_x(t < V_{E}) \le e^{-\rho t}/2.
\end{equation*}
By induction on $k\ge 1$ with the Markov property,
we deduce that for any $k$ :
\begin{equation*}
	\frlq{x\in (0, 1)}
	\PR_x(k\, t < V_{E}) \le e^{-k\rho t}/2^k.
\end{equation*}
We know conclude the proof of Lemma \ref{TS:lem:esc}
by noting :
\begin{align*}
		\E_x(\exp[\rho V_{E}])
		&\le \Tsum{k\ge 0} e^{\rho t\, [k+1]} \PR_x(V_{E} \in [k\, t, (k+1)\,t) )
		\\&\le e^{\rho t} \Tsum{k\ge 0} 2^{-k} = 2 e^{\rho t} <\infty.
\end{align*}
\epf
\subsection{Proof of Proposition \ref{TS:p:r1r0ra} :
	convergence to $\delta_1$ for $\rho_1<\rho_0<\rho_\alpha$ }

For $\mu \in \MnxTS$, 
we have the following
lower-bound of the mass absorbed at $1$ before time 1 :
\begin{equation*}
	\mu P_1 \{1\} 
	= \PR_\mu( \tau_1 \le 1\le \extO )
	\ge \xi\, \PR_{1/n} (\tau_1 \le 1)\, \exp[-\Ninf{r}]
\end{equation*}
By the Markov property, 
this implies with 
$C_{n,\xi}
:= \xi\, \PR_{1/n} (\tau_1  \le 1)\, \exp[- (\Ninf{r} - \rho_1)]$~:
\begin{equation*}
	\mu P_t \{1\} 
	= \mu P_1 \{1\} \exp[-\rho_1 (t-1)]
	\ge C_{n,\xi} \, \exp[-\rho_1\, t].
	\EQn{TS:Cnx1}{}
\end{equation*}
Since $\tau_0 \le \extOU$ and the extinction rate is 
$\rho_0$ once $X$ has reached $0$, 
then exploiting \Req{TS:etaB}~: 
\begin{align*}
		&\PR_\mu (\tau_0 \le t < \ext) 
		\le \E_\mu\lc
		\exp[-\rho_0\,(t - \extOU)]\pv
		\extOU \le t\rc
		\\& \hcm{1}
		\le  \exp[-\rho_0\,t]\, 
		\lc 1 + \rho_0\, \intO{t}ds\, \exp[\rho_0\, s] \PR_\mu(s<\extOU)  \rc
		\\&\hcm{1}
		\le \exp[-\rho_0\,t]\, 
		\lc 1 + \Ninf{\bar{\heig}}\times \frac{\rho_0}{\ra - \rho_0}\rc.
		\EQn{TS:majr0}{ }
\end{align*}
With again \Req{TS:etaB}, and \Req{TS:Cnx1}
and the fact that both $\mu A_t$ and $\delta_1$ are probability measure :
\begin{align*}
		&\NTV{\mu A_t - \delta_1}
		= \sup_{D} |\mu A_t(D) - \delta_1(D)|
		= \mu A_t[0, 1)
		= \mu P_t[0,1) / (\mu P_t[0,1) + \mu P_t \{1\})
		\\&\hcm{1}
		\le \dfrac{\PR_\mu (\tau_0 \le t < \ext)  + \PR_\mu (t < \extOU) }
		{\mu P_t \{1\}}
		\le C'_{n, \xi}\, \exp[-(\rho_0-\rho_1)\,t]
		\\&
		\where C'_{n, \xi} 
		:=  \lc  1+  \Ninf{\bar{\heig}}\times \ra / (\ra - \rho_0)\rc 
		/  C_{n,\xi}
		\SQ
\end{align*}

\subsection{Proof of Proposition \ref{TS:p2:r1r0ra} : 
	conditional convergence to $\delta_0$ for $\rho_0<\rho_\alpha$}
Let $t\ge 1$ and assume first 
that $\mu([0, x]) \ge \xi$ for $x\in (0,1)$ and $\xi>0$.
From \Req{TS:etaB},
\begin{align*}
		&
		\PR_\mu (t < \extOU) 
		\le \Ninf{\bar{\heig}}\, \exp[-\ra\, t]
\end{align*}
With the rough lower-bound 
$\mu\, P^1_1\{0\} \ge \exp[-\Ninf{r}]\, \PR_\mu (\tau_0 \le 1)$~: 
\begin{align*}
		&\hcm{1}
		\PR_\mu (\tau_0 \le t < \extU) 
		\ge \exp[-\Ninf{r}]\, \PR_\mu (\tau_0 \le 1)\times \exp[-\rho_0\, (t-1)]
		\\&\hcm{3}
		\with C
		:= \dfrac{\Ninf{\bar{\heig}}\,\exp[\Ninf{r} - \rho_0]} 
		{\xi\, \PR_x (\tau_0 \le 1)} >0,
		\\&
		\NTV{\mu A^{1}_t
			- \delta_0} 
		= \mu A^{1}_t(0,1)
		= \dfrac{\PR_\mu( t < \extOU) }
		{\PR_\mu( t < \extOU) + \PR_\mu (\tau_0 \le t < \extU) }
		\\&\hcm{3}
		\le C\,  \exp[-(\ra-\rho_0)\, t].
		\EQn{TS:ra-r0}{}
\end{align*}

The case where $\mu$ has support on $\{0,1\}$ is trivial, 
since then 
$\mu A^1_t = \delta_0$.

Finally, for the general case of $\mu \in \MOne \setminus \{\delta_1\}$, 
where $\mu(0,1) >0$, 
remark that, for any $s>0$,  
there exists $m_s \in (0,1)$ such that ~: 
\begin{align*}
		&\mu A^1_s = m_s\, \mu A^{01}_s + (1-m_s)\, \delta_0
		\\&\text{ where for any } x>0
		\mVg \qquad
		\mu A^{01}_s([0, x])
		\cvifty{s}{ \alpha([0, x]) >0}.
\end{align*}
by Proposition \ref{TS:p:QSD01}, 
with the rate of convergence  uniform over $\mu$.
Thus, 
we deduce some $t_\vee>0$ such that, with $x = 1/2$ ~:
\begin{align*}
		\frlq{\mu \in \MOne \setminus \{\delta_1\}}
		\mu\, A^{1}_{t_\vee}([0, x])
		\ge \mu A^{01}_{t_\vee}([0, x])
		\ge \alpha([0, x]) /2 := \xi.
\end{align*}
Thus, for any $t\ge t_\vee$, 
by the Markov property, then \Req{TS:ra-r0} 
with the initial condition $\mu A^{01}_{t_\vee}$ :
\begin{align*}
		\NTV{\mu A^{1}_t
			- \delta_0} 
		= \NTV{\mu A^{1}_{t_\vee}\, A^1_{t-{t_\vee}}
			- \delta_0} 
		\le C\, \exp[(\ra-\rho_0)\,t_\vee] \exp[-(\ra-\rho_0)\, t].
		\SQ
\end{align*}

\subsection{Proof of Proposition \ref{TS:p:r1rar0} and \ref{TS:p2:r1rar0}~:
	the case $\rho_1<\rho_\alpha<\rho_0$}
For this Proposition, we need to adapt 
the proof given in Subsection 5.3
of \cite{AV_QSD}.
The main step is to prove that the mass on the interval $(0,1)$ does not vanish~:
\begin{lem}
	\label{TS:l:mxt}
	Assume that $\ra < \rho_0$. Then, there exists $n_{\circ}\ge 2, \xi_{\circ} >0$ such that~:
	\begin{align*}
			&\hspace{1cm}
			\frl{n\in \N}
			\frl{\xi >0}
			\Ex{t_{\circ} > 0}
			\\&\hspace{.5cm}
			\frl{\mu \in \MnxTS}
			\frl{t\ge t_{\circ}}\qquad 
			\mu A^1_t(1/n_{\circ},\,1-1/n_{\circ}) \ge \xi_{\circ}.
	\end{align*}
\end{lem}

We also need to ensure the persistence of the component issued from the coupling,
which is done with the following lemma.
\begin{lem}
	\label{TS:l:ps1}
	Assume that $\ra < \rho_0$ and $\zeta \in \MoneM$. 
	Then, there exists $t_{P}, c_{P}>0$ such that~:
	\begin{align*}
			\frl{x\in [0,1)}
			\frlq{t\ge t_{P}}
			\PR_x(t<\extU) 
			\le c_{P}\, \PR_{\zeta}(t<\extU).
			\EQn{TS:cps}{}
	\end{align*} 
\end{lem}

The measure $\zeta$ comes from a mixing estimate that we recall 
--cf \Req{TS:alc} :
\begin{lem}
	\label{TS:l:mix1}
	Let $n\ge 2$ and $\xi >0$. 
	Then, there exists $\zeta \in \MoneM$, 
	$t_{M}, c_{M}>0$ such that
	for any $\mu \in \MOne$ satisfying $\mu(1/n,\,1-1/n) \ge \xi$~:
	\begin{align*}
			\mu A^1_{t_{M}}(dx) 
			\ge c_{M}\, \zeta(dx).
	\end{align*}
\end{lem}

\subsubsection{Proof of Proposition \ref{TS:p:r1rar0}}
The proof of Proposition \ref{TS:p:r1rar0} follows as 
the one of Proposition \ref{TS:p2:r1r0ra} from \Req{TS:etaB1},
which is a consequence of Proposition \ref{TS:p2:r1rar0}.
\epf

\subsubsection{Proof of Proposition  \ref{TS:p2:r1rar0} with Lemmas \ref{TS:l:mxt}-3}

Combining these tree lemmas and applying exactly the same reasoning 
as in Subsection 5.3 of \cite{AV_QSD}
proves that there exists a unique QSD $\alpha_1$ associated to $\extU$,
with the convergence stated in Proposition  \ref{TS:p2:r1rar0}.

Moreover, 
as stated in Proposition  \ref{TS:p2:r1rar0},
we can identify $\alpha_1$ and $\heig^1$.

Let $\alpha^y_1 := y\, \alpha + (1-y)\, \delta_0$. 
For any $t\ge 0$~:
\begin{align*}
		&\alpha^y_1 P^1_t(dx)
		= y\, \exp[-\ra\, t]\, \alpha(dx)
		+  \lc (1-y)\,\exp[-\rho_0\, t] 
		+ y\, \PR_\alpha( \tau_0 \le t< \extU)\rc
		\delta_0(dx).
		\EQn{TS:aly1}{}
\end{align*}
It follows from Theorem 2.6 in \cite{coll}
that
the exit state is independent from the exit time 
when the initial condition is a QSD, 
with an exponential law for the exit time.
Thus~:
\begin{align*}
		& \PR_\alpha( \tau_0 \le t< \extU)
		=\E_\alpha \lc \exp[-\rho_0\, (t - \tau_0)] \pv
		\tau_0 = \extOU \le t\rc 
		\\&\hcm{2}
		= \PR_{\alpha}\lp \tau_0 = \extOU\rp\,
		\E_\alpha \lc \exp[-\rho_0\, (t - \extOU)] \pv
		\extOU \le t\rc
		\\&\hcm{2}
		= \PR_{\alpha}\lp \tau_0 = \extOU\rp\, 
		\intO{t} \exp[-\rho_0\, (t - s)]\times \ra\, \exp[-\ra\, s]\, ds
		\\&\hcm{2}
		= \lp \exp[- \ra\, t] - \exp[-\rho_0\, t]\rp
		\dfrac{\ra\,\PR_{\alpha}\lp \tau_0 = \extOU\rp}{\rho_0 - \ra}.
		\EQn{TS:rho0al}{}
\end{align*}
With our choice \Req{TS:y0}, i.e.
$\frac{1- y_{\alpha}}{y_\alpha} = \frac{\ra\,\PR_{\alpha}\lp \tau_0 = \extOU\rp}{\rho_0 - \ra}$,
we see that we obtain indeed : 
\begin{align*}
		\alpha^{y_\alpha}_1 P^1_t = \exp[-\ra\, t] \, \alpha^{y_\alpha}_1.
\end{align*}

The proof that $\heig^1$ is uniquely defined, 
the convergences in \Req{TS:alU} and \Req{TS:CVeta1}
and the upper-bound in \Req{TS:etaB1} are exactly the same as in \cite{AV_QSD}.
It remains to identify $\heig^1$. 
Clearly~:
\begin{equation*}
	\heig^1(0) := \limInf{t}  \exp[\rho_\alpha\, t]\, \PR_0(t<\extU) 
	=  \limInf{t}  \exp[-(\rho_0 - \ra)\, t] = 0.
	\EQn{TS:eta10}{}
\end{equation*}
We define $\heig^1_t$ similarly as $h_t$ (cf \Req{TS:ht}) :
for $x \in [0, 1)$
\begin{equation*}
	\heig^1_t(x) :=  \exp[\rho_\alpha\, t]\, \PR_x(t<\extU).
	\EQn{TS:eta1}{\heig^1_t}
\end{equation*}
Decomposing according to the state of $X$ at time $1$, 
and recalling that since $X$ stays at $0$ 
once it is hit,
we have :
\begin{equation*}
	\heig^1_{2\,t}(x) 
	= \heig_{t}(x)\, 
	\LAg \delta_x\, A^{01}_t\bv
	\heig^1_t\RAg
	+ \heig^1_t(x)\, \delta_x\, A^1_t\{0\}
	\,\heig^1_t(0)
	\EQn{TS:eta2t}{}
\end{equation*}
From \Req{TS:eta10} and \Req{TS:etaB1}, 
the second term in the right-hand side is clearly negligible.
Because of \Req{TS:CVal}, 
we are interested in the asymptotic as $t\ifty$
of $\LAg \alpha\bv \heig^1_t\RAg$.
Since $\alpha$ is the QSD associated to the extinction rate $\ra$
with extinction at time $\extOU$, 
we already know that $\PR_\alpha (t< \extOU) = \exp[-\ra t]$.
We will deal with the component that has fixed at $0$ before time $t$
thanks to \Req{TS:rho0al}.
Concluding with \Req{TS:y0},
we mean that~:
\begin{align*}
		&\LAg \alpha\bv
		\heig^1_t\RAg
		= \exp[\rho_\alpha \, t] \; \PR_{\alpha} (t< \extOU)
		+  \exp[\rho_\alpha \, t] \;\PR_{\alpha}( \tau_0 \le t< \extU) 
		\\&\hcm{2}
		= 1 + \lp 1- \exp[-(\rho_0-\ra)\, t]\rp
		\dfrac{\ra\,\PR_{\alpha}\lp \tau_0 = \extOU\rp}{\rho_0 - \ra}
		\\&\hcm{3}
		\cvifty{t}
		{1 + (1 - y_\alpha) / y_\alpha = 1/ y_\alpha}. 
		\EQn{TS:alEta1}{}
\end{align*}
From \Req{TS:eta2t},  \Req{TS:alEta1}, \Req{TS:CVal},   and \Req{TS:etaB1},
we conclude~:
\begin{align*}
		\heig^1_{2\,t}(x) 
		\cvifty{t}{\heig(x) / y_\alpha 
			= \heig^1(x)}.
		\SQ
\end{align*}


\subsubsection{Proof of Lemma \ref{TS:l:ps1}}
From \Req{TS:CVeta1}, with the notation \Req{TS:eta1} : 
\begin{equation*}
	\PR_{\zeta}(t< \extU) 
	= \LAg \zeta\bv \heig^1_t\RAg \; \exp[-\ra\, t]
	\where 
	\LAg \zeta\bv \heig^1_t\RAg
	\cvifty{t}{\LAg \zeta \bv \heig\RAg / y_\alpha}.
\end{equation*}

We know from Proposition \ref{TS:p:QSD01}
that $h$ is lower-bounded by a positive constant on any $\cD_n$
for $n \ge 2$. 
It implies in particular $\LAg \zeta \bv \heig\RAg>0$.
From \Req{TS:CVal}, let thus $t_{P} >0$ be such that :
\begin{equation*}
	\frlq{t\ge t_{P}}
	\LAg \zeta\bv \heig^1_t\RAg 
	\ge \LAg \zeta \bv \heig\RAg / 2 > 0. 
	\EQn{TS:minal1}{t_{P}}
\end{equation*}
\Req{TS:eta10} is clearly true and implies with \Req{TS:etaB}\, 
that \Req{TS:cps} holds for $x =0$.

For $x\in (0,1)$
and any $t>0$~:
\begin{align*}
		& \PR_x( \tau_0 \le t< \extU)
		=\E_x \lc \exp[-\rho_0\, (t - \tau_0)] \pv
		\tau_0 = \extOU \le t\rc 
		\\&\hcm{1}
		\le \E_x \lc \exp[-\rho_0\, (t - \extOU)] \pv
		\extOU \le t\rc 
		\\&\hcm{1}
		=  \exp[- \rho_0\, t] \, 
		\lp 1+ \rho_0\,\intO{t} \exp[\rho_0\, s]
		\times \PR_x(s\le \extOU \le t)\, ds
		\rp
		\\&\hcm{1}
		\le \exp[- \rho_0\, t] \, 
		\lp 1+ \rho_0\,\Ninf{\bar{\heig}}\,
		\intO{t} \exp[(\rho_0-\ra)\, s]  ds \rp
		\\&\hcm{1}
		\le  \exp[- \rho_0\, t] \,  
		+ \dfrac{\rho_0\,\Ninf{\bar{\heig}}\,}{\rho_0-\ra}\,  \exp[- \ra\, t].
		\EQn{TS:maj0}{}
\end{align*}
Combining \Req{TS:minal1}, \Req{TS:maj0} and \Req{TS:etaB}\, ends the proof
of Lemma \ref{TS:l:ps1}.
\epf

\subsubsection{Proof of Lemma \ref{TS:l:mxt} }
Let $n_{\circ}\ge 3$ such that~:
\begin{align*}
		\alpha(1/n_{\circ}\mVg 1- 1/n_{\circ} ) \ge 1/2
		\EQn{TS:nxt}{n_{\circ}}
\end{align*}
From \Req{TS:CVal}, we can find $t_{S} >0$ such that
for any $\mu$ with $\mu(0,1)>0$~:
\begin{align*}
		&\hcm{0.5}
		\frlq{t\ge t_{S}}
		\PR_\mu (t< \extOU) 
		\ge \LAg \mu\bv \heig \RAg /2 \times \exp[-\ra\, t]
		\mVg
		\\&
		\NTV{\mu A^{01}_t -\alpha} \le 1/4
		\quad \text{ thus  } \quad
		\mu A^{01}_t(1/n_{\circ}\mVg 1- 1/n_{\circ} ) \ge 1/4.
		\EQn{TS:tsb}{t_{S}}
\end{align*}

Since $0$ is absorbing and by $\Req{TS:tsb}$ :
\begin{align*}
		&
		\frlq{t\ge 0}
		\mu A^1_t (dx)
		= \mu A^1_t(0,1)  \times \mu A^{01}_t(dx)
		+ \lc 1- \mu A^1_t(0,1) \rc\, \delta_0(dx),
		\\&\hcm{.5}  
		\frlq{t\ge t_{S}}
		\mu A^{1}_t (1/n_{\circ}\mVg 1- 1/n_{\circ} )
		\ge \mu A^1_t(0,1) / 4
		\where :
		\EQn{TS:tnX}{}
		\\&
		\mu A^1_t(0,1) 
		= \dfrac{\PR_\mu(t<\extOU)}{\PR_\mu(t<\extOU)+ \PR_\mu(\tau_0 \le t<\extU)}
		= \lp 1 
		+ \dfrac{\PR_\mu(\tau_0 \le t<\extU)}{\PR_\mu(t<\extOU)} 
		\rp^{-1}.
		\EQn{TS:muA1}{}
\end{align*}

Assume first that $\mu[1/n, 1-1/n] \ge \xi$ for some $n\ge 3$ and $\xi >0$.
Since $\heig$ is positive on $(0,1)$,
this implies, with \Req{TS:CVal}, \Req{TS:maj0}, \Req{TS:tnX}\, and \Req{TS:muA1}, 
a lower-bound
$\xi_{\circ}$ that only depends on $n$ and $\xi$ such that~:
\begin{align*}
		\frlq{t\ge t_{S}}
		\mu A^{1}_t (1/n_{\circ}\mVg 1- 1/n_{\circ} )
		\ge \xi_{\circ}.
		\EQn{TS:nxth}{}
\end{align*}

Lemma \ref{TS:esc1}
completes the proof.
Indeed, consider  $\mu \in \MnxTS$ 
(w.l.o.g. $\mu\{1\} = 0$ since it vanishes immediately).

Either $\mu(1/n, x_\vee) \ge \xi /2$ 
and we deduce the result from \Req{TS:nxth},

or $\mu[x_\vee, 1) \ge \xi /2$ 
and we deduce from Lemma \ref{TS:esc1} and \Req{TS:nxth}~:
\begin{align*}
		\frlq{t\ge t_{S}+t_{B}}
		\mu A^{1}_t (1/n_{\circ}\mVg 1- 1/n_{\circ} ) 
		= [\mu A^{1}_{t_{B}}] A^{1}_{t-t_{B}} (1/n_{\circ}\mVg 1- 1/n_{\circ} ) 
		\ge \xi'_{\circ}.
		\SQ
\end{align*}

\subsection{Proof of Proposition \ref{TS:p:r1r0a}~:
	the case $\rho_1<\rho_0 = \rho_\alpha$}

The calculations leading to \Req{TS:maj0} gives for the case $\rho_0 = \ra$~:
\begin{align*}
		\frlq{\mu \in \MOne}
		\PR_\mu( \tau_0 \le t< \extU)
		\le \exp[- \rho_0\, t] \, 
		\lp 1+ \rho_0\,\Ninf{\bar{\heig}}\,t \rp.
		\EQn{TS:maj1}{}
\end{align*}
With \Req{TS:maj1} instead of \Req{TS:majr0},
like in the proof of Proposition \ref{TS:p:r1r0ra}
(i.e. with \Req{TS:etaB}\, and \Req{TS:Cnx1}),
we deduce Proposition \ref{TS:p:r1r0a}.
\epf

\subsection{Proof of Proposition \ref{TS:p2:r1r0a}~:
	conditional convergence to $\delta_0$ 
	when $\rho_0 = \rho_\alpha$}

Thanks to Proposition \ref{TS:p:QSD01},
there is for any $n\ge 2$ a positive lower-bound of $\heig$ in $\cD_n$,
$\LAg \mu \bv \heig\RAg$ is uniformly lower-bounded 
for $\mu \in \MnxTS^{0,1}$ (for any $n\ge3,\xi>0$).
By \Req{TS:CVal}, 
for $t$ sufficiently large 
and any $\mu \in \MnxTS^{01}$~:
\begin{align*}
		\PR_\mu (t<\extOU) 
		\ge c_{n,\xi}\, \exp[- \rho_0\, t].
\end{align*}
Combining this with \Req{TS:muA1} and \Req{TS:maj1} concludes the proof
that for $t$ sufficiently large ~:
\begin{align*}
		\NTV{\mu A^1_t - \delta_0} 
		\ge C_{n,\xi} / t.
		\SQ
\end{align*}
\textbf{Remark : } 
Adapting the proof of step 1, Section 5.1 of \cite{Ch1D},
one can prove that there exists $t_B, c_B, n_B>0$ such that 
for any $x\in [1/2, 1)$, $\delta_x A^1_{t_{B}}\in [1/4, 1-1/n_B)$.
Exploiting the above proof of Proposition \ref{TS:p2:r1r0a}, 
this implies that the convergence is uniform 
for any $\mu$ such that $\mu [1/n, 1) \ge \xi$.
\\

For the reverse inequality,
assume first that $\mu \in \MnxTS^{01}$.
From the Markov property
\begin{align*}
		&\PR_\mu( \tau_0 \le t< \extU)
		= \E_\mu \lc
		\exp[- \rho_0\, (t- \tau_0)] \, \pv
		\tau_0 = \extOU \le t\rc
		\\&\hcm{1}
		= \E_\mu \lc
		\exp[- \rho_0\, t] \, 
		\lp 1+ \rho_0\,\intO{\extOU} \exp[\rho_0\, s]\, ds\rp\pv
		\tau_0 = \extOU \le t\rc
		\\&\hcm{1}
		\ge \rho_0\, 
		\exp[- \rho_0\, t] \,
		\int_{0}^{t} \exp[\rho_0\, s]\, 
		\PR_\mu(\tau_0 = \extOU \in [s, t]) \, ds,
		\\& \hcm{1}
		\ge c_{n, \xi} 
		\exp[- \rho_0\, t] \,
		\int_{0}^{t} 
		\PR_{\mu A^{01}_s} (\tau_0 = \extOU \le t-s)\, ds,
		\EQn{TS:tau0ext}{}
\end{align*}
where we exploited once more Proposition \ref{TS:p:QSD01}
in the last inequality,
to obtain a uniform lower-bound $c_{n, \xi}$ 
on $\LAg \mu \bv \heig_s\RAg$.

Since $\PR_\alpha (\tau_0 = \extOU)>0$ and
by monotone convergence,
there exists $t_\vee>0$ such that~:
\begin{align*}
		\frlq{t\ge t_\vee}
		\PR_\alpha (\tau_0 = \extOU \le t) 
		\ge \PR_\alpha (\tau_0 = \extOU \le t_\vee) 
		:= m_0 >0.
		\EQn{TS:tvee}{t_\vee}
\end{align*}
Now, according to \Req{TS:CVal},
we choose $t_{S}>0$ such that~:
\begin{align*}
		& \frl{\mu \in \MoneM}
		\frlq{s\ge t_{S}}
		\NTV{\mu A^{01}_s - \alpha} 
		\le m_0 / 2
		\\
		\text{which implies } \quad 
		&\frl{s\ge t_{S}}\frlq{t-s\ge t_\vee}
		\PR_{\mu A^{01}_s} (\tau_0 = \extOU \le t-s) 
		\ge m_0/2.
		\EQn{TS:tsb2}{t_{S}}
\end{align*}
Thus, \Req{TS:tsb2} and \Req{TS:tau0ext} imply 
that for any $t\ge t_{S} + t_\vee$~:
\begin{align*}
		&\PR_\mu( \tau_0 \le t< \extU)
		\ge c'_{n, \xi}\, \exp(-\rho_0\, t)\times (t - t_{S} - t_\vee).
\end{align*}
With \Req{TS:etaB} and \Req{TS:muA1}, this concludes the proof 
that~:
\begin{align*}
		\mu[1/n, 1- 1/n] \ge \xi
		\quad \imp\quad
		\frlq{t\ge  t_{S} + t_\vee}
		\NTV{\mu A^1_t - \delta_0} 
		\le C_{n,\xi}/t.
		\EQn{TS:cSt}{}
\end{align*}

Now, we prove that such upper-bound is in fact uniform 
with respect to $\MoneM$
thanks to Lemma \ref{TS:esc1}.
Indeed 
\begin{align*}
		&\mu A^1_{t_{B}} (dx)
		= \mu A^1_{t_{B}}(0,1)  \times \mu A^{01}_{t_{B}}(dx)
		+ \lc 1- \mu A^1_{t_{B}}(0,1) \rc\, \delta_0(dx),
		\\&\where 
		\Ex{\xi_{B} >0}\Exq{n_{B} \ge 2}
		\frlq{\mu \in \MoneM}
		\mu A^{01}_{t_{B}}(1/n_{B}, 1- 1/n_{B}) \ge \xi_{B}
		\\&
		\text{Thus by \Req{TS:cSt}~: } \hcm{1}
		\frlq{t\ge t_{S}+t_\vee}
		\NTV{[\mu A^{01}_{t_{B}}] A^1_t  - \delta_0 }
		\le  c_{n_{B}, \xi_{B}}/t.
		\EQn{TS:NTVce}{}
\end{align*}
We also note that there exists $y_t\in (0,1)$ such that~:
\begin{align*}
		&\mu A^1_{t_{B}+t} (dx)
		= y_t \, [\mu A^{01}_{t_{B}}] A^1_t 
		+ (1- y_t) \delta_0.
\end{align*}
In fact, our comparison of the survival from $0$ and from $\mu$
gives us a uniform upper-bound $C>0$ such that~:
\begin{align*}
		& y_t = \mu A^1_t(0,1) \times 
		\dfrac{\LAg \mu A^{01}_{t_{B}} \bv
			\PR_.(t<\extU)\RAg}
		{\LAg \mu A^{1}_{t_{B}} \bv
			\PR_.(t<\extU)\RAg}
		\quad
		\le C\, \mu A^1_t(0,1)
\end{align*}
Hence, we have more precision on the convergence :
\begin{equation*}
	\NTV{\mu A^{1}_{t +t_{B}} - \delta_0 }
	\le \mu A^1_{t_{B}}(0,1) \times C / t.
\end{equation*}

And at least, \Req{TS:NTVce} concludes the proof of Proposition \ref{TS:p2:r1r0a} 
(where $t_{S} +t_\vee$ replaces $t_\vee$).
\epf
\\

\subsection{Proof of Proposition \ref{TS:p:r10ra}~:
	the case $\rho_0 = \rho_1<\rho_\alpha$}
Since $\rho_0 = \rho_1$, it is straightforward that
any convex combination of  $\delta_0$ and $\delta_1$ 
is a QSD, with extinction rate  $\rho_1$.

It is then not difficult to adapt the proof of Proposition \ref{TS:p2:r1r0ra},
and since $\PR_\mu(\tau_{0,1} \le 1)$ is lower-bounded
uniformly over any $\mu \in \MOne$, we obtain
\begin{align*}
		\frlq{\mu \in \MOne}
		\mu A_t (0,1) 
		\le C\,  \exp[-(\ra-\rho_0)\, t].
\end{align*}
\begin{align*}
		&\mu A_t\{0\}
		= \dfrac{\E_\mu\lc \exp[-\rho_1\, (t-\extOU)]\pv 
			\extOU = \tau_0\le t\rc}
		{\PR_\mu(t< \extOU)
			+\E_\mu\lc \exp(-\rho_1\, (t-\extOU))\pv 
			\extOU = \tau_{0,1}\le t\rc}
		\\&
		= \dfrac{\E_\mu\lc \exp[\rho_1\, \extOU]\pv 
			\extOU = \tau_0\le t\rc}
		{\E_\mu\lc \exp[\rho_1\, \extOU]\pv 
			\extOU = \tau_{0,1}\le t\rc}
		\times \lp 1 + \dfrac{\exp[\rho_1\, t]\, \PR_\mu(t< \extOU)}
		{\E_\mu\lc \exp(\rho_1\, \extOU)\pv 
			\extOU = \tau_{0,1}\le t\rc} \rp^{-1}.
		\EQn{TS:muA0}{}
\end{align*}
The limit as $t \ifty$ is well-defined 
and the convergence occurs at exponential rate since~:
\begin{align*}
		&0\le 
		\E_\mu\lc \exp[\rho_1\, \extOU]\pv 
		\extOU = \tau_1\rc
		- \E_\mu\lc \exp[\rho_1\, \extOU]\pv 
		\extOU = \tau_1\le t\rc
		\\&\hcm{1}
		\le \E_\mu\lc \exp[\rho_1\, \extOU]\pv 
		t <\extOU \rc
		\\&\hcm{1}
		\le \Ninf{\bar{\heig}}\, \exp[-\ra\, t]\,
		\lc 1 + \rho_1 \intRp \exp[\rho_1\, s]\, \PR_{\mu A^{01}_t} (s < \extOU)\, ds \rc
		\\&\hcm{1}
		\le \Ninf{\bar{\heig}}\, 
		\lc 1 + \frac{\rho_1\, \Ninf{\bar{\heig}}}{\ra - \rho_1}\rc\;
		\exp[-\ra\, t] 
		:= C\,  \exp[-\ra\, t].
\end{align*}
The same holds of course for the case $\Lbr \extOU = \tau_{0,1}\Rbr$
and $\E_\mu\lc \exp(\rho_1\, \extOU)\pv 
\extOU = \tau_{0,1}\le t\rc$ 
converges with exponential rate.
Therefore with \Req{TS:muA0} --and the well-defined notation \Req{TS:xmu}--
we can define some $C>0$ such that $\forall \mu \in \MOne$~:
\begin{align*}
		|\mu A_t\{1\} - (1-x(\mu))|
		\vee  |\mu A_t\{0\} - x(\mu)|
		\vee |\mu A_t(0,1)|
		\le C\, \exp[-(\ra-\rho_0)\, t],
\end{align*}
which concludes the proof of Proposition \ref{TS:p:r10ra}.
\epf

\subsection{Proof of Proposition \ref{TS:p:rar10}~:
	the case $\ra < \rho_0\wedge \rho_1$}
This proof is very similar to the one of Proposition \ref{TS:p2:r1rar0},
so we won't go into much detail.
Lemmas \ref{TS:l:mix1} and  \ref{TS:l:ps1}
are of course replaced by~:
\begin{lem}
	\label{TS:l:mxt2}
	Assume that $\ra < \rho := \rho_0\wedge \rho_1$. 
	Then, there exists $n_{\circ}\ge 3, \xi_{\circ} >0$ such that~:
	\begin{align*}
			&\hspace{1cm}
			\frl{n\in \N}
			\frl{\xi >0}
			\Ex{t_{\circ} > 0}
			\\&\hspace{.5cm}
			\frl{\mu \in \MnxTS^{01}}
			\frl{t\ge t_{\circ}}\qquad 
			\mu A_t(1/n_{\circ},\,1-1/n_{\circ}) \ge \xi_{\circ}.
	\end{align*}
\end{lem}

\begin{lem}
	\label{TS:l:ps2}
	Assume that $\ra < \rho := \rho_0\wedge \rho_1$ and $\zeta \in \MoneM$. \\
	Then, there exists $t_{P}, c_{P}>0$ such that~:
	\begin{align*}
			\frl{x\in [0,1]}
			\frlq{t\ge t_{P}}
			\PR_x(t<\ext) 
			\le c_{P}\, \PR_{\zeta}(t<\ext).
	\end{align*} 
\end{lem}

We leave the proofs to the reader, 
and just mention that we can take as an upper-bound for
$\PR_x( \tau_1 \le t< \ext)$ 
the same formula as for $\PR_x( \tau_0 \le t< \ext) = \PR_x( \tau_0 \le t< \extU)$,
with $\rho_1$ instead of $\rho_0$ 
(cf \Req{TS:maj0}).

For the rest of the proof, 
we remark that,
for $\alpha^y := y_\alpha\, \alpha + y_0\, \delta_0 + y_1\, \delta_1$
with $y_\alpha +y_0+y_1 = 1$,
\Req{TS:aly1} has to be replaced by~:
\begin{align*}
		&\alpha^y P_t(dx)
		= y_\alpha\, \exp[-\ra\, t]\, \alpha(dx)
		+  \lc y_0\,\exp[-\rho_0\, t] 
		+ y_\alpha\, \PR_\alpha( \tau_0 \le t< \extU)\rc
		\delta_0(dx)
		\\&\hcm{2}
		+ \lc y_1\,\exp[-\rho_1\, t] 
		+ y_\alpha\, \PR_\alpha( \tau_1 \le t< \extU)\rc
		\delta_1(dx).
		\EQn{TS:aly01}{}
\end{align*}
Again : $\alpha^y P_t(dx) = \exp[-\ra\, t]\; \alpha^y(dx)$ 
iff the conditions in \Req{TS:a01} are satisfied.
\epf

\subsection{Proof of Proposition \ref{TS:p:r1ar0}~:
	the case $\ra = \rho_1 < \rho_0$}
Let us first prove that we only need to control 
$\NTV{\mu A^0_t  - \delta_1}$ 
like it is done in Proposition \ref{TS:p2:r1r0a}.
From Proposition \ref{TS:p2:r1rar0},
we know that 
for some $\alpha_1 := y_\alpha\, \alpha + y_0\, \delta_0$, 
with $y_\alpha, y_0 \in (0,1)$, 
there exists $C^1, \chi^1 >0$ such that~:
\begin{align*}
		\NTV{\mu A^1_t - \alpha_1} 
		\le C^1\, \exp[-\chi^1\, t].
		\EQn{TS:zeta1}{\chi_1}
\end{align*}

Consequently, for $t$ sufficiently large~:
\begin{align*}
		\dfrac{ y_0}{2 \,y_\alpha}
		\le \dfrac{ \mu A_t \{0\} }{\mu A_t(0,1)}
		\le \dfrac{2 \, y_0}{y_\alpha}.
		\EQn{TS:0vsAl}{}
\end{align*}

On the other hand,
with the notation 
$\quad
\mu A^0_t (dx)
:= \PR_\mu \lp X_t \in dx	\bv 	t< \extO\rp$~:
\begin{align*}
		& \NTV{\mu A_t - \delta_1}
		= \lc 1 + \dfrac{\mu A_t \{1\}} 
		{\mu A_t (0,1) +  \mu A_t \{0\}} \rc^{-1}
		\mVg \hcm{1}
		\NTV{\mu A^0_t - \delta_1}
		= \lc 1 + \dfrac{\mu A_t \{1\}} {\mu A_t (0,1)} \rc^{-1}.
\end{align*}
Consequently, 
\Req{TS:0vsAl} implies that $\NTV{\mu A_t - \delta_1}$
has the same rate of convergence 
as $\NTV{\mu A^0_t  - \delta_1}$ 
(as long as it indeed converges to $0$).

Now, from the proof of Proposition \ref{TS:p2:r1r0a},
we deduce quite immediately~:
\begin{align*}
		&\hcm{1.2}
		\Exq{t_\vee, C>0}
		\frl{t\ge t_\vee}
		\frlq{\mu \in \MOne}
		\NTV{\mu A^0_t
			- \delta_1}
		\le C  / t,
		\\& \frl{n\ge 3}
		\frlq{\xi>0}
		\Ex{t_{n,\xi}, c_{n, \xi}>0}
		\frl{t\ge t_{n, \xi}}
		\frlq{\mu \in \MnxTS^{0,1}}
		\\&\hcm{3}
		\NTV{\mu A^0_t
			- \delta_1}
		\ge c_{n, \xi}  / t.
		\SQ
\end{align*}

\subsection{Proof of Proposition \ref{TS:p:r10a}~:
	the most critical case $\ra = \rho_0 = \rho_1$}
Any convex combination of $\delta_0$ and $\delta_1$ 
is clearly a QSD with extinction rate $\rho := \rho_0 = \rho_1 = \rho_\alpha$.

%
For $t \ge 0$ and $x\in [0,1]$, let :
\begin{align*}
		& \heig_t (x) := \exp[\rho\, t] \, \PR_x(t < \extOU)
		\EQn{TS:EtaT}{\heig_t}
		\mVg
		E^t_0(x) :=  \E_x\lc \exp[\rho\, \extOU]\pv 
		\extOU = \tau_0\le t\rc
		\quad
		\\&  E^t_1(x) :=  \E_x\lc \exp[\rho\, \extOU]\pv 
		\extOU = \tau_1\le t\rc,
		\EQn{TS:Et1}{E^t_1}
\end{align*}

Let then $k\ge 1$ and $\mu \in \MOne$ with $\mu(0,1) >0$,
so that $\LAg \mu\bv \heig \RAg >0$. Then~:
\begin{align*}
		&\LAg \mu \bv E^{k}_0 \RAg
		= \sum_{j=0}^{k-1}
		\LAg\mu \bv \heig_{j}\RAg\; 
		\LAg \mu A^{01}_{j} \bv E^1_0\RAg,
\end{align*}
where by \Req{TS:CVal}, 
(with the upper-bound $e^\rho$  of $E^1_0$),
there exists $C >0$ such that~:
\begin{align*}
		|\LAg\mu \bv \heig_{j}- \heig\RAg\;| 
		\le C\, \exp[-j\, \chi]
		\mVg \quad
		|\LAg \mu A^{01}_{j} - \alpha \bv E^1_0\RAg|
		\le C \,  \exp[-j\, \chi].
\end{align*}
Consequently~:
\begin{align*}
		&|\LAg \mu \bv E^{k}_0 \RAg
		- k\, \LAg\mu \bv \heig\RAg\; 
		\LAg \alpha \bv E^1_0\RAg|
		\le 2\, C / (1+ \exp[-\chi]) 
		< \infty.
		\EQn{TS:muEk}{}
		\\\text{ Likewise }
		&|\LAg \mu \bv E^{k}_0 + E^{k}_1 \RAg
		- k\, \LAg\mu \bv \heig\RAg\; 
		\LAg \alpha \bv E^{1}_0 + E^{1}_1\RAg|
		\le 4\, C / (1+ \exp[-\chi]) 
		< \infty.
\end{align*}
From \Req{TS:muA0} and \Req{TS:etaB}, 
we deduce that there exists $C'>0$ such that~:
\begin{align*}
		&\left| \mu A_k\{0\} 
		- \dfrac{\LAg \alpha \bv E^{1}_0\RAg}
		{\LAg \alpha \bv E^{1}_0 + E^{1}_1\RAg}\right|
		\le \dfrac{C'}{k\, \LAg\mu \bv \heig\RAg}
		\EQn{TS:muaK}{}
\end{align*}
The symmetrical result for $\mu A_k\{1\}$ 
holds of course true, 
and since the sum of the limits equals 1,
we deduce also 
\begin{align*}
		|\mu A_k(0,1)| \le \dfrac{C'}{k\, \LAg\mu \bv \heig\RAg}.
		\EQn{TS:muAk}{}
\end{align*}

Again, from Theorem 2.6 in \cite{coll},
the exit state is independent from the exit time 
when the initial condition is a QSD, 
with an exponential law for the exit time. 
Thus~:
\begin{align*}
		\LAg \alpha \bv E^{1}_0\RAg
		= \PR_\alpha \lp \tau_0 = \extOU\rp
		\; \intO{1} \exp[\rho\, s]\, \rho\, \exp[-\rho\, s]\, ds
		= \PR_\alpha \lp \tau_0 = \extOU\rp.
\end{align*}
To end the proof, 
just remark that $\LAg\mu \bv \heig\RAg$ is lower-bounded for any $\mu \in \MnxTS^{01}$.
\epf

\subsection{Proof of Proposition \ref{TS:prop:SI}~:
	$\rho_\alpha(\gamma) \rightarrow \infty$ as 
	$\gamma \rightarrow \infty$}

We assume first that $r\equiv 0$
and choose arbitrary some $t$, for instance $t := 1$.
Consider $T_\delta := \inf\{u\ge 0\pv X_u\, (1-X_u) \le \delta\}$,
which can possibly be 0.
Given any $\eps >0$,
we want to prove that choosing $\delta$ sufficiently small
ensures, uniformly for $\gamma \ge 1$ :
$\PR_x(T_\delta \le t_0 \mVg 2\, t_0 < \tau_{0,1}) \le \eps$.

We can notice that~:
\begin{align*}
		X(t_0) = x_0 - s\,.\, T(t_0) + \gamma\, \wtd{B}[T(t_0)]
		\with  T(t_0) := \int_0^{t_0\wedge \tau_{0,1}} X_u\, (1-X_u)\, du,
		\EQn{TS:Tch}{}
\end{align*}
and $\wtd{B}$ has the law of a Brownian Motion.
Indeed, define 
$$\wtd{B}_v := B(T^{-1}(v)\wedge \tau_{0,1})
+\idc{\tau_{0,1}< T^{-1}(v)}\, 
(\hat{B}(T^{-1}(v)) - \hat{B}(\tau_{0,1})),$$
with $\hat{B}$ another Brownian Motion independent of $B$.
Since for any $v>0$, $T^{-1}(v) := \inf\{t\ge 0 : 
\int_0^t X_u\, (1-X_u)\, du > v\}$
is a stopping times,
$\hat{B}$ is indeed a continuous martingale 
with respect to the filtration $\F_{T^{-1}(v)}$.
Finally, 
the change of variable $w = T(u)$ ensures that 
$\E[ (\wtd{B}_v - \wtd{B}_{v'})^2] = v - v'$ for any $v>v'$.

%
%

On the other hand, by Itô's formula :
\begin{align*}
		&\E_x \lp X_t - x - \int_0^{t_0\wedge T_\delta} s\, X_u\, (1-X_u)\, du\rp^2
		= \E_x \lp \int_0^{t_0\wedge T_\delta} \gamma\, \sqrt{ X_u\, (1-X_u)}\, dB_u\rp^2
		\\&= \E_x \lc \int_0^{t_0\wedge T_\delta} \gamma^2\, X_u\, (1-X_u)\, du\rc
		\\ &\text{thus }
		\PR_x(t_0< T_\delta)\times t_0\, \gamma^2\, \delta 
		\le (2+s\, t_0/2)^2, \quad \text{independent from $x$.}
\end{align*}
For $\gamma$ sufficiently large,
it implies that
$\PR_x(t_0< T_\delta)$ is indeed lower than $\eps$.

Thus $\PR_x( t_0 < \tau_{0,1}) 
\le \PR_x(T_\delta \le t_0 \mVg 2\, t_0 < \tau_{0,1})
+ \PR_x(t_0< T_\delta)\le 2\eps$.

In the general case of bounded $r$, 
we deduce for the QSD $\alpha$ that,
for $\gamma$ large enough :
\begin{align*}
		\PR_{\alpha}( 2\,t_0 < \tau_{0,1}) 
		= \exp[-2\,\rho_\alpha\, t_0]
		\le 2\, \exp[2\,\|r\|_\infty\, t_0]\, \eps.
\end{align*}
It indeed proves that $\rho_\alpha \rightarrow \infty$ 
as $\gamma \rightarrow \infty$.
\epf

\subsection{Proof of Proposition \ref{TS:prop:rI}~:
	$\rho_\alpha < \rho_0\wedge \rho_1$ for $r$ sufficiently strong}

Define $r_2, r_3$ such that $\max r(x) < r_3 < r_2 <r(1)\wedge r(0)$
and the open sets $A := r^{-1}([0, r_3))\subset B := r^{-1}([0, r_2)) \subset (0,1)$
(recall that $r$ is assumed to be continuous).
We choose arbitrary $t_0$. 
A classical result on diffusion ensures that there exists $\rho>0$ such that :
\begin{align*}
		\inf_{x\in A} \PR_x(X_{t_0} \in A\mVg \frl{s\le t} X_s \in B)
		\ge \exp[-\rho\, t_0].
\end{align*}
Then, it implies by the Markov property~:
\begin{align*}
		\inf_{x\in A} \PR_x(\frl{s\le t} X_s \in B\pv t_0 < \extOU)
		\ge C \, \exp[-(\rho+ R\, r_2)\, t].
\end{align*}
From the Harnack inequality, 
we know that $\alpha^{(R)}$ 
has a lower-bounded density on any open set of $(0,1)$
so that $\alpha^{(R)}(A) >0$ and 
\begin{align*}
		\alpha^{(R)}(B)
		\ge \exp[\rho_\alpha^{(R)}\,t] \PR_{\alpha^{(R)}} (\frl{s\le t} X_s \in B\pv t < \extOU)
		\ge C \,\alpha^{(R)}(A)  \exp[-(\rho+ R\, r_2- \rho_\alpha^{(R)})\, t].
\end{align*}
This proves $\rho_\alpha^{(R)} \le \rho + R\, r_2 
< R\,( r(0)\wedge r(1)) = \rho_0^{(R)}\wedge \rho_1^{(R)}$
for $R$ sufficiently large.
\epf

\subsection{Proof of Proposition \ref{TS:prop:r0}~:
	$ \rho_0\wedge \rho_1<\rho_\alpha $ for $r$ sufficiently weak}

Let $\rho_\alpha^{(0)}$ be the death rate of the QSD 
for the Wright-Fisher diffusion conditioned not to touch the boundary
with $r = 0$.
Since $\rho_0^{(R)}\wedge \rho_1^{(R)} 
= R\times [r^0(0)\vee r^0(1)] \cvz{R} 0$,
it is sufficient to prove that $\liminf_R \rho_\alpha^{(R)} \ge \rho_\alpha^{(0)}$
to deduce Proposition \ref{TS:prop:r0}.

By Proposition \ref{TS:p:QSD01}, for any $x\in (0,1)$ and $t\ge 0$ :
\begin{align*}
		&\PR_x(t<\tau_{0,1}) \le \Ninf{\bar{\heig}}\, \exp[-\rho_\alpha^{(0)}\, t]
		\PR_x(t<\tau_{\partial, 0,1}) \le \Ninf{\bar{\heig}} \, \exp[-(\rho_\alpha^{(0)}-R\,\|r\|_\infty)\, t]
\end{align*}
and in particular, with the QSD  $\alpha(R)$ as initial condition,
we deduce $\rho_\alpha^{(R)} \ge \rho_\alpha^{(0)}-R\,\|r\|_\infty$.

\epf

\textsl{Remark : }
In fact, $\rho_\alpha^{(R)} \cvz{R}{\rho_\alpha^{(0)} >0}$, because :
\begin{align*}
		\exp[-(\rho_\alpha^{(0)}+R\,\|r\|_\infty - \rho_\alpha^{(R)})\, t]
		\le \exp[+(\rho_\alpha^{(R)})\, t],
		\PR_{\alpha^{(0)}}(t<\tau_{\partial, 0,1}) 
		\cvifty{t}{\LAg \alpha^{(0)}\bv \heig^{(R)}\RAg},
\end{align*}
which implies  $\rho_\alpha^{(R)} \le \rho_\alpha^{(0)} + R\,\|r\|_\infty$.

\subsection{Proof of Proposition \ref{TS:prop:sI}~:
	concentration towards $0$ as $\gamma \rightarrow 0$}

Since $r$ is bounded, 
the probability of the event $\Lbr t < \tau_\eps\Rbr$ with $r$ 
is at most $\exp(\|r\|_\infty\, t)$ times 
the probability with $r \equiv 0$.
If we prove that the latter converges to 0 
(as a limit of this parameter $\gamma$),
it will be the same for the former.
We can thus assume without loss of generality that $r \equiv 0$.

We recall (see \Req{TS:Tch}) that for any $t\ge 0$
and initial condition $1 - \eps$~:
\begin{align*}
		X(t) = 1- \eps - s\,.\, T(t) + \gamma\, \wtd{B}[T(t)]
		\with  T(t_0) := \int_0^{t_0\wedge \tau_{0,1}} X_u\, (1-X_u)\, du,
\end{align*}
and $\wtd{B}$ has the law of a Brownian Motion.

Fix some $M>0$
ans assume that $\gamma \le \eps / (2\,M)$ 
and that we are conditionally on the event 
$\{\sup_{u\le 1/s} |\wtd{B}_u| \le M\}$.
Then, $T(\infty)< 1/s$, since $T^{-1}(1/s)$ would be well-defined otherwise and would satisfy :
\begin{equation*}
	X(T^{-1}(1/s) ) \le 1- \eps - 1 + \eps/2 \le -\eps/2 <0,
\end{equation*}
which contradicts the classical property 
that $X$ takes its value on $[0, 1]$.
By the definition of $T(t)$,
we see that there exists $c = c(\eps)>0$ such that~:
\begin{equation}\EQn{TS:tinf}{}
	1/s > T(\tau_{\eps/2}\wedge \tau_{1-\eps/2})
	\ge c\, \tau_{\eps/2}\wedge \tau_{1-\eps/2}.
\end{equation}
We also deduce that for any $u\le T(\infty)$,
since $u< 1/s$~:
\begin{align*}
		X(T^{-1}(u)) \le 1- \eps  + \eps/2
		\le 1-\eps/2.
\end{align*}
This implies that $\tau_{1-\eps/2} = \infty$ on the event $\{\sup_{u\le 1/s} |\wtd{B}_u| \le M\}$
and by \Req{TS:tinf} that $\tau_{\eps/2} \le t := 1/(c s)$.
This directly implies~:
\begin{align*}
		\frlq{\gamma \le  \eps / (2\,M)}
		\PR_{1-\eps}(t < \tau_\eps) 
		\le \PR (\textstyle \sup_{u\le 1/s} |B_u| \ge M).
\end{align*}
Letting $M$ tend to $\infty$ concludes the proof of Proposition \ref{TS:prop:sI}.
\epf

\section{Proof of Theorem \ref{TS:MLim}}
\label{TS:PMLim}
\setcounter{eq}{0}

Like in \cite{LM15}, 
the proof follows a standard procedure 
\cite{FM04}, \cite{JM86}, \cite{CFM06} 
in which we prove : (i) the tightness of the sequence of stochastic processes
-- which implies a subsequential limit, 
and (ii) the uniqueness of this limit. 
For the tightness of $\{\mu_{t}^{m,n}\}_{m,n}$
on $D([0,\ T],\ \mathcal{P}([0,1]))$,
it is sufficient, by Theorem 14.26 in Kallenberg \cite{K97}
to show that $\{\LAg \mu_{t}^{m,n}\bv f\RAg\}$ is tight on
$D([0,\ T],\ \R)$ for any test function $f$ 
from a countably dense subset of continuous, positive
functions on $[0$, 1$]$.

\subsection{Semimartingale property of multilevel selection process}

It will be useful for
what follows to treat $\LAg \mu_{t}^{m,n} \bv f\RAg $ as a semimartingale. 
We exploit the following discrete derivatives of $f$,
with span $1/n$ :
\begin{align*}
		&D_{x}^{+}f(x) := n\, (f(x+1/n)  - f(x))\mVg\quad
		D_{x}^{-}f(x) := n\, (f(x)  - f(x-1/n))\mVg
		\\&
		D_{xx}f(x) := n^2 \ltm (f(x+1/n) + f(x-1/n) - 2 f(x)) 
		= n \ltm (D_{x}^{+}f(x) - D_{x}^{-}f(x)).
\end{align*}

We recall that in our limit, 
$n, m \rightarrow \infty$,
$\bar{\gamma}_I/n \rightarrow \gamma_I$,
$\bar{\gamma}_I\, \bar{s} \rightarrow s$,
$\{\bar{\gamma}_G \bar{r}(x)\} $ is fixed
and  $\bar{\gamma}_G$ bounded.
It is easy to adapt the proof of  \cite{LM15} in order to state~: 
\begin{lem}
	For $f\in C^{2}([0,1])$ and $\mu_{t}^{m,n}$  with generator $L^{m,n}$  defined in (1),
	\begin{align*}
			\LAg  \mu_{t}^{m,n}\bv f\RAg 
			-\LAg  \mu_{0}^{m,n} \bv f \RAg \ 
			=A_{t}^{m,n}(f) 
			+M_{t}^{m,n}(f)
	\end{align*}
	where $A_{t}^{m,n}(f)$  is a process of finite variation, $A_{t}^{m,n}(f) :=\int_{0}^{t}a_{s}^{m,n}(f)ds$, with~:
	\begin{align*}
			&a_{t}^{m,n}(f)
			=  \sum_{i}
			\mu_{t}^{m,n}\,  \lp ^i\!/_n\rp
			\times\ 
			^i\!/_n
			\lp 1- ^i\!/_n \rp\ 
			\lc  \lp ^{\bar{\gamma}_I}\!/_{n}\rp\, D_{xx}f\lp ^i\!/_n\rp 
			-  \bar{\gamma}_I\, \bar{s} \ D_{x}^{-}f\lp ^i\!/_n\rp \rc
			\\&\hcm{1}
			+\bar{\gamma}_G\ \Lbr\sum_{j}\mu_{t}^{m,n}\lp ^j\!/_n\rp
			\bar{r}\lp ^j\!/_n\rp \, f\lp ^j\!/_n\rp
			-\sum_{i} \mu_{t}^{m,n} \lp ^i\!/_n\rp 
			f\lp ^i\!/_n\rp 
			\sum_{j}\mu_{t}^{m,n}
			\lp ^j\!/_n\rp   \bar{r}\lp ^j\!/_n\rp \Rbr
	\end{align*}
	and $M_{t}^{m,n}(f)$ is a c\`{a}dl\`{a}g martingale with (conditional)  quadratic variation~:
	\begin{align*}
			&  \LAg M^{m,n}(f)\RAg _{t}
			=\frac{1}{m}\int_{0}^{t}
			\Lbr\frac{\bar{\gamma}_I}{n} \sum_{i} 
			\mu_{s}^{m,n} \lp ^i\!/_n\rp
			\ ^i\!/_n\ 
			\lp 1-^i\!/_n \rp\ 
			\lc	\lp 
			D_{x}^{+}f \lp ^i\!/_n\rp 		\rp^{2}
			+(1+\bar{s})\lp D_{x}^{-}f	\lp ^i\!/_n\rp 
			\rp^{2}\rc
			\right.
			\\&\hcm{2}
			\left.
			+\bar{\gamma}_G\ \sum_{i,j} \mu_{s}^{m,n} \lp ^i\!/_n\rp
			\mu_{s}^{m,n} \lp ^j\!/_n\rp\  \lp 1+\bar{r}(^j\!/_n)\rp
			\lc f\lp ^i\!/_n\rp  
			- f\lp ^j\!/_n\rp \rc^{2}\Rbr 
			ds
	\end{align*}
\end{lem}

\subsection{Proof of the convergence to our limit}

We prove here that the drift term is tight 
while the martingale converges to zero.



For the finite variation term $A_{t}^{m,n}(f)$,
assuming w.l.o.g. $\bar{\gamma}_I/n \le 2\, \gamma_I$,
$\bar{\gamma}_I\, \bar{s}\le 2\, s$ :
\begin{align*}
		&|a_{t}^{m,n}(f)|\ 
		\leq\ 2\gamma_I \Ninf {f''}
		+ 2\, s\ \Ninf {f'}
		+2\ \Ninf {r} \Ninf{f}
		:= G_f
		\\&\hcm{3}
		\text{ therefore : }
		\sup_{t\in[0,T]}|A_{t}^{m,n}(f)|\ \leq G_{f}T
\end{align*}
where $G_{f}$ is a constant that depends on $f$. 
Moreover,
for any prescribed $\eps$,
we can always choose $\delta_\vee$ 
to be sufficiently small 
so that, 
for any $0\le t\le t +\delta$ 
with $\delta\le \delta_\vee$,
for any $n, m$~:
$|A_{t+\delta}^{m,n}-A_{t}^{m,n}| 
\le \delta G_{f}
\leq\eps.$ 
By Proposition 3.26, Chapter 3 in \cite{JS03}, 
this proves immediately 
that the sequence $(A_{t}^{m,n})_{t\le T}$ is tight,
and any limit is continuous.
\\
For the martingale part, 
assuming w.l.o.g. $\bar{s}\le 1$
~:
\begin{align*}
		\LAg M_{t}^{m,n}(f)\RAg_{t}
		\leq\ \frac{T}{m}\Lbr 
		6\, \gamma_I\, \Ninf {f'}
		+ (\bar{\gamma}_G+ \Ninf{r})\,
		\Ninf{f}^2 \Rbr 
		:= J_{f}/m \cvifty{m}{0},
\end{align*}
where $J_f$ is a constant 
only depending on $T >0$ and $f\in C^{2}([0,1])$. 
From Burkholder-Davis-Gundy’s inequality,
since the jumps of $M_{t}^{m,n}(f)$ 
are bounded by $\Ninf{f}/m$ :
\begin{align*}
		\E\lc \sup\!_{t\le T} (M_{t}^{m,n}(f))^2\rc
		\le C\, J_{f}/m + \Ninf{f}^2/m^2 
		\cvifty{n, m}{0}.
\end{align*}
This proves that $M_{t}^{m,n}(f)$ converges to $0$
in such a way that 
$(\LAg  \mu_{t}^{m,n}\bv f\RAg 
-\LAg  \mu_{0}^{m,n} \bv f \RAg)_{t\le T}$
is tight 
and any associated limit is continuous.

By construction and the Law of Large Numbers, $\LAg  \mu_{0}^{m,n} \bv f \RAg$ 
converges to $\LAg  \mu_{0} \bv f \RAg$.
Thus,  the sequence $(\mu_{t}^{m,n})_{t\le T}$ for $n, m \ge 1$
is tight in $D([0; T]; \M_1([0;1]))$.

So, we consider a subsequence $(\mu^{(k)}_t)_{t\le T} = (\mu^{m_k, n_k}_t)_{t\le T}$
such that $m_k, n_k \ifty$ as $k\ifty$
and such that $(\mu^{(k)}_t)$ converges to $(\mu_t)_{t\le T}$
in $D([0; T]; \M_1([0;1]))$.
Necessarily, $\mu_0$ coincide with the law of the initial condition
provided in the assumptions of Theorem \ref{TS:MLim}.
For any $f \in C_2([0,1])$ and $t\le T$, 
it is not difficult to see that as $k\ifty$~:
\begin{align*}
		&\LAg  \mu^{(k)}_t\bv f\RAg \rightarrow \LAg  \mu_t\bv f\RAg
		\mVg \quad
		\LAg  \mu^{(k)}_0\bv f\RAg \rightarrow \LAg  \mu_0\bv f\RAg
		\mVg
		\\&
		a_{t}^{(k)}(f) 
		\rightarrow  \LAg  \mu_t\bv \cL_{W\!F} f\RAg
		+ \LAg  \mu_t\bv r \ltm f\RAg - \LAg  \mu_t\bv f \RAg\ltm \LAg  \mu_t\bv r\RAg
		\mVg \quad 
		M_{t}^{(k)}(f) \rightarrow 0.
\end{align*}
Thus $(\mu_t)$ is a solution to equation \Req{TS:eqCar}.
From the uniqueness property that we proved in Proposition \ref{TS:Char},
and the tightness of the sequence,
we conclude that $(\mu^{m, n}_t)_{t\le T}$
converges globally to this solution.
This concludes the proof of Theorem \ref{TS:MLim} 
and more globally the proofs presented in this paper.
\epf
\\

\textbf{Acknowledgment}\\
This paper would not have been the same without the support 
of Etienne Pardoux, my PhD supervisor. I sincerely wish to address him my thanks
for this article as well.

	\typeout{get arXiv to do 4 passes: Label(s) may have changed. Rerun}

\begin{thebibliography}{00}
	\addcontentsline{toc}{chapter}{Bibliographie}
	
	\bibitem[CFM06]
	{CFM06}
	Champagnat, N., Ferri\`ere, R., M\'{e}l\'{e}ard, S.;
	Unifying evolutionary dynamics: from individual stochastic processes to macroscopic models. Theoretical Population Biology, Elsevier, V.69, N.3, pp.297-321 (2006)
	
	\bibitem[CH19]
	{Henry} 
	Champagnat, N., Henry, B.; 
	A probabilistic approach to Dirac concentration
	in non-local models of adaptation with several resources,
	Ann. Appl. Probab.,
	V.29, N.4, pp.2175-2216 (2019)
	
	
	\bibitem[CV18a]
	{Ch1D} Champagnat, N., Villemonais, D.;
	Uniform convergence of conditional distributions for absorbed one-dimensional diffusions, 
	Adv. in Appl. Probab., V.50, N.01, pp.178-203 (2018)
	
	\bibitem[CV17b]
	{ChpLyap} Champagnat, N.,  Villemonais, D.;
	Lyapunov criteria for uniform convergence of conditional distributions of absorbed Markov processes,
	preprint on ArXiv~: 1704.01928
	(2017)
	
	\bibitem[CMS13]
	{coll}  Collet, P., Martínez, S., San Martín, J.; 
	Quasi-Stationary Distributions, 
	Probab. and Its Appl., Springer, Berlin Heidelberg (2013)
	
	\bibitem[Daw10]
	{D10} 
	Dawson, D.A.; 
	Introductory Lectures on Stochastic Population Systems,
	Technical Report Series No. 451, Laboratory for Research in Stat. and Probab. (2010)
	
	
	\bibitem[DZ98]
	{DZ98}
	Dembo, A., Zeitouni, O.; Large Deviation Techniques and Applications. 
	Second Edition. Springer-Verlag, New York  (1998).
	
	\bibitem[Ew04]
	{Ew04} Ewens, W.; 
	Mathematical Population Genetics, I, 
	2nd ed., Interdiscip. Appl. Math. 27, Springer-Verlag, Berlin (2004)
	
	\bibitem[FM04]
	{FM04}
	Fournier, N., Méléard, S.; A microscopic probabilistic description of a locally regulated
	population and macroscopic approximations. Ann. Appl. Probab., V.14, N.4, pp.1880-1919 (2004)
	
	\bibitem[Gr06]{G06} 
	Grafen, A.; Optimisation of inclusive fitness;
	J. Theor.  Biol., V.238, pp.541–563 (2006)
	
	\bibitem[JS03]
	{JS03} Jacod, J., N.Shiryaev, A.; 
	Limit Theorems for Stochastic Processes.
	Second Edition. Grundlehren der mathematischen Wissenschaft, V.288, Springer (2003)
	
	\bibitem[JS17]{JS17} 
	Jenkins, P. A., Spano, D:  Exact simulation of the Wright-Fisher diffusion. 
	Ann. Appl. Prob., V.27, N.3, pp.1478-1509 (2017)
	
	
	
	\bibitem[JM86]{JM86}
	Joffe, A., Metivier, M.; 
	Weak convergence of sequences of semimartingales with applications to multitype branching processes. 
	Advances in Appl. Probab., V.18, pp.20-65 (1986)
	
	
	\bibitem[Ka97]
	{K97}  Kallenberg, O.; Foundations of Modern Probability. Springer, 1st edition (1997) 
	
	\bibitem[LK$^+$07]{L+07}
	Lehmann, L., Keller, L., West, S. A., Roze, D.; Group selection
	and kin selection. Two concepts but one process. 
	Proceed. of the  National Acad. of Sciences of the U.S.A., V.104,
	pp.6736–6739 (2007)
	
	\bibitem[Luo13]{L13}
	Luo, S.; A unifying framework reveals key properties of multilevel selection, J. Theor.
	Biol. V.341, pp.41-52 (2013)
	
	\bibitem[LM15]{LM15} 
	Luo, S., Mattingly, J.; 
	Scaling limits of a model for selection at two scales. 
	Nonlinearity. V.30, N.4, pp.1682-1707 (2015)
	
	\bibitem[ML14]{ML14} 
	Mullon, C., Lehmann, L.; 
	The robustness of the weak selection approximation for the evolution
	of altruism against strong selection. 
	J. Evol. Biol., V.27, pp.2272-2282
	(2014) 
	
	\bibitem[O'F08]{O08} O'Fallon, B.;
	Population Structure, Levels of Selection, and the Evolution of Intracellular Symbionts; 
	Evolution,
	V.62, N.2, pp.361-373 (2008)
	
	\bibitem[RG$^+12$]{RGHK12}
	Roques, L., Garnier, J., Hamel, F., Klein, E.K.;
	Allee effect promotes diversity in traveling waves of colonization,
	PNAS June 5,  V.109, N.23, pp.8828-8833 (2012)
	
	\bibitem[SGE13]
	{SGE13} Schraiber, J., Griffiths, R. C., Evans, S. N.;
	Analysis and rejection sampling of Wright-Fisher diffusion bridges.
	Theoretical Population Biology, V.89, pp.64–74 (2013)
	
	
	\bibitem[TN06]{TN06} Traulsen, A., Nowak, M.A.: 
	Evolution of cooperation by multilevel selection.
	Proceedings of the National Academy of Sciences, 103, 29 (2006)
	
	\bibitem[Ve19]{AV_QSD} Velleret, A.; 
	Unique Quasi-Stationary Distribution, with a stabilizing extinction, preprint available on ArXiv at :
	\url{https://arxiv.org/abs/1802.02409}
	
	\bibitem[WD19]{WD19} Wang, Z., Durrett, R. J.;
	Extrapolating Weak Selection in Evolutionary Games,
	Math. Biol.  V.78, N.1-2, pp.135-154 (2019)
	
	
	\bibitem[WGG07]{W+07} West, S.A., Griffin, A.S., Gardner, A.; 
	Social semantics: altruism, cooperation,  mutualism, strong reciprocity and group selection. 
	J. Evol. Biol. V.20,  pp.415–432 (2007)
	
	\bibitem[WGG07b]{W+08} West, S.A., Griffin, A.S., Gardner, A.;  
	Social semantics: how useful has group selection been? 
	J. Evol. Biol. V.21, pp.374–385 (2007)
	
	\bibitem[WMG11]
	{W+11}
	West, S.A., El Mouden, C., Gardner, A.;
	Sixteen common misconceptions about the evolution of cooperation in humans, 
	Evol. and Human Behavior V.32, N.4, pp.231-262 (2011)
	
	\bibitem[WG$^+$13]{WGHT13}
	Wu, B., Garcia, J., Hauert, C., Traulsen, A.;
	Extrapolating Weak Selection in Evolutionary Games, 
	PLoS Comput. Biol.,  V.9, N.12, pp.1-7 (2013) 
	
\end{thebibliography}
\end{document}